\documentclass{article}
\usepackage{tikz,pdfsync}
\usepackage{xcolor}
\usepackage[english]{babel}
\usepackage{textcomp, marvosym, amsfonts, amsmath, amssymb}
\usepackage{pdfsync,stmaryrd}
\usepackage{times}
\usepackage{graphicx}

\def\scbrl{[\hspace{-.14em}[}
\def\scbrr{]\hspace{-.14em}]}
\def\scbr#1{\scbrl#1\scbrr}

\newcommand{\gknf}{=_{k\mbox{-nf }}}
\newcommand{\pk}{\psi_k}
\newcommand{\xk}{\xi_k}

\newcommand{\pskpe}{\psi_{k+1}}
\newcommand{\xskpe}{\xi_{k+1}}

\newcommand{\ck}{\chi_k}
\newcommand{\cskpe}{\chi_{k+1}}

\newcommand{\Om}{{\Omega}}

\newcommand{\al}{{\alpha}}

\newcommand{\be}{{\beta}}

\newcommand{\ga}{\gamma}

\newcommand{\de}{\delta}
\newcommand{\eo}{{\varepsilon_0}}

\newcommand{\la}{\lambda}

\newcommand{\om}{\omega}

\def\scbrl{[\hspace{-.14em}[}
\def\scbrr{]\hspace{-.14em}]}
\def\scbr#1{\scbrl#1\scbrr}

\newtheorem{theo}{Theorem}
\newtheorem{lem}{Lemma}
\newtheorem{defi}{Definition}

\providecommand{\keywords}[1]{\textbf{{Keywords---}} #1}

\begin{document}
\title{Giant and illusionary giant Goodstein principles}

\author{Andreas Weiermann\\ Ghent University\\ Department of Mathematics: Analysis, Logic and Discrete Mathematics\\ Krijgslaan 281 S8\\ 9000 Ghent, Belgium\\
{\tt email: Andreas.Weiermann@UGent.be}}

\date{}
\maketitle
\begin{abstract}
We analyze several natural Goodstein principles which themselves are defined with respect to the Ackermann function and the extended Ackermann function.
These Ackermann functions are well established canonical fast growing functions labeled by 
ordinals not exceeding $\varepsilon_0$. 
Among the Goodsteinprinciples under consideration, the giant ones, will be proof-theoretically strong (being unprovable in $\mathrm{PA}$ in the Ackermannian case
and being unprovable in $\mathrm{ID}_1$ in the extended Ackermannian case)
whereas others, the illusionary giant ones, will turn out to be comparatively much much weaker although they look strong
at first sight.
\end{abstract}
\section{Introduction}
This article is part of a general program on exhibiting natural independence results for first order Peano arithmetic $\mathrm{PA}$, its subsystems and 
strong extensions, as for example, the theory $\mathrm{ID}_1$ of non iterated inductive definitions. (For readers more familiar with reverse mathematics let us remark that he theory $\mathrm{ID}_1$ has the same proof-theoretic strength as $\mathrm{ACA}_0+(\Pi^1_1-\mathrm{CA})^-$
where $(\Pi^1_1-\mathrm{CA})^-$ refers to the scheme of lightface $\Pi^1_1$-comprehension.)
We are mainly interested in natural combinatorial statements which have non trivial unprovability strength. 
The Goodstein principles serve as canonical examples since they provide prime examples of the intrinsic interplay between finite and infinite numbers.

Goodstein's original result \cite{Goodsteina,Goodsteinb} dealt with iterated exponential terms and their canonical interpretation in the ordinals less than $\eo$.
The basic idea is as follows. There are two processes involved. One is making a given number under normal circumstances greater
and the other making a given number smaller by subtracting a one. Before round $0$ a non negative integer is chosen.
Then like in a two person game the following moves are carried out alternatingly. 
At round $k\geq 0$ the first player develops the number, which is already obtained in the game, completely with respect to base $k+2$ as long as it is not yet zero.
Here also the exponents are developed hereditarily. Then he replaces every occurrence of $k+2$ by $k+3$. 
Then the second player makes an apparently innocent move.
She just subtracts a one from the previous result. Then the game moves over to round $k+1$ with the new result.
The second player wins this specific game if the number zero is reached after finitely
many rounds. 
Surprisingly the second player always wins but this fact is unprovable in $\mathrm{PA}$. (See, for example, \cite{Kirby,Cichon} for a proof.)
Goodstein's theorem is an example for a giant Goodstein principle 
since the underlying sequences become very (unprovably) long even for small starting values.

In this paper we deal with Goodstein principles which are defined with respect to 
canonical representations of natural numbers using the Ackermann function and the extended Ackermann function.
The paper therefore splits naturally into two parts. 

In the fist part we deal with normal forms and Goodstein principles which are defined relative to the Ackermann function.
 Similar principles have already been studied in \cite{Weiermann} and \cite{Arai} for the Ackermann function.
 These articles are based on a very complicated iterated sandwiching procedure and it is quite natural to ask what happens
 if the sandwiching is reduced to a one step approximation. This question has been investigated for the Ackermann function which starts at the bottom level with the exponential 
 function in \cite{Fernandez}. It turned out that the strength of the resulting Goodstein principles dropped considerably and we arrived in these cases
 at intermediate Goodstein principles.
 
In the first part we investigate normal forms which are based on a one step approximation with respect to the Ackermann function
which starts at the bottom level with the successor function. Surprisingly the resulting Goodstein principles is still not provable from the 
axioms of first order arithmetic when base changes are carried out
in the two critical arguments of the Ackermann function which is involved in the number representations. 
This principle leads thus to a strong independence result and we call it a giant Goodstein principle.
This reflects the fact that it will match in strengths with the strongest possible Goodstein principle which is based on number representations defined in terms of the Ackermann function.
In \cite{walk} it has been shown that the latter principle is equivalent to the one consistency of first order arithmetic $\mathrm{PA}$.

But we also consider restrictions where the base change is carried out only with respect to one critical argument 
of the Ackermann function used to build up normal form representations under consideration. These lead to Goodstein principles which look strong at first sight
but as a closer analysis reveals they are in fact of low proof-theoretic strength. One principle will be equivalent over primitive recursive arithmetic with the assertion
that $\om^{\om+\om}$ does not admit primitive recursive descending sequences. 
The other will be equivalent with the one consistency of primitive recursive arithmetic.

We call the latter principle an  illusionary giant Goodstein principle.\footnote{These principles show similarities with a prominent figure from the literature, namely Tur Tur, known from the story: Jim Button and Luke the Engine Driver written by Michael Ende.
"After a long and hazardous journey, Jim Button and Luke the Engine Driver arrive in the Dragon City. Along the way, they make two new friends, 
the giant Mr. Tur Tur (who is actually a "Scheinriese", an "illusionary giant"  he only appears to be a giant from afar; when approached, it turns out he is actually of normal height), and Nepomuk, the half-dragon." [Citation modelled after https://enacademic.com/dic.nsf/enwiki/2981127]}

In the second technically much more advanced part we consider Goodstein principles defined with respect to the extended Ackermann functions which are canonical fast growing functions indexed by
 ordinals less than $\eo$. 

It turns out that the resulting Goodstein principle becomes very strong when
base changes are performed in two critical arguments of the the extended Ackermann function. 
In fact a termination proof necessarily involves a detour via uncountable ordinals since 
it becomes unprovable in $\mathrm{ID}_1$. For proving this we develop a novel theory of majorization properties of fundamental sequences.
This theory allows us to prove that $k$-normal forms after base change are moved into $k+1$-normal forms.

Moreover we employ the machinery of Buchholz's collapsing function $\psi$. 
We found it very amazing that our ordinal mapping from natural numbers into ordinals has the property that normal forms for numbers 
are moved by magic to ordinals in Buchholz $\psi$ normal form and we believe that this underpins the naturality of our approach.

As far as we know this example will be the first example of a Goodstein principle for natural numbers which has such a high unprovability strength.

In analogy with the first part one would expect that the strength of the Goodstein principle will drop seriously when base change is defined with respect to one critical argument of the extended Ackermann function used to define the number representations.
Quite surprisingly it turned out that a base change only in the first critical argument still suffices to generate a giant Goodstein principle.

But a base change in the second critical argument only leads to an illusionary Goodstein principle which in strength is equivalent with the classical Goodstein principle.

\section{Giant and illusionary giant Goodstein sequences related to the Ackermann function}

Let us define the Ackermann function $A_a(k,b)$ with respect to iteration parameter $k<\om$ as follows:

\begin{eqnarray*}
A_0(k,b)&:=&b+1,\\
A_{a+1}(k,0)&:=&A_a(k,\cdot)^k(0),\\
A_{a+1}(k,b+1)&:=&A_a(k,\cdot)^k(A_{a+1}(k,b)).
\end{eqnarray*}

Here the upper index denotes the number of function iterations.
It is a routine matter to show that for any fixed $k\geq 1$ the function $a,b\mapsto A_a(k,b)$ is not primitive recursive whereas for
fixed $a$ the function $b\mapsto A_a(k,b)$ is primitive recursive. 
It is also easy to show that the function $k,a,b\mapsto A_a(k,b)$ is strictly monotone in $a,b$ and $k\geq 1$.

{\bf Convention.} From now on $k$ denotes, if not stated otherwise, a positive integer not smaller than $3$.
\begin{lem}

For all $m>0$ there exist unique $a,b,l<\om$ such that
\begin{enumerate}
\item $m=A_a(k,b)+l,$
\item $a$ is maximal with $A_a(0)\leq m$ (so that $A_a(k,0)\leq m <A_{a+1}(k,0)$),
\item $b$ is maximal with $A_a(b)\leq m$ (so that $A_a(k,b)\leq m <A_{a}(k,b+1)$).
\end{enumerate}
\end{lem}
We write $m\gknf A_a(k,b)+l$ in this case.
If $m\gknf A_a(k,b)+l$ and $a=0$ then necessarily $m=A_0(k,\cdot)^m(0)$ where $m<k$. 
If $m=k$ then $m\gknf A_1(k,0)$. Indeed we see by an easy induction on b that $A_1(k,b)=k\cdot (1+b)$.
This yields in particular that $A_a(k,b)> 2\cdot b$ for $a>0$.
The latter estimate is used tacitly at several occasions.

In the sequel we often write $A_a(b)$ for $A_a(k,b)$ and $B_a(b)$ for $A_a(k+1,b)$ when $k$ is fixed in a given context.

\begin{lem}
\begin{enumerate}
\item If $m\gknf A_a(b)+l$ and $l>0$ then $A_a(b)+l-1$ is in $k$ normal form, too.
\item If $m\gknf A_a(b)$ and $b>0$ then  $A_a(b-1)$ is in $k$ normal form, too.
\item If $m\gknf A_a(0)$ then  for $0<l<k$ $A_a^l(0)$ is in $k$ normal form, too.
\end{enumerate}
\end{lem}
Proof. All assertions are easy to see. Let us prove the last assertion. If  $0<l<k$ then $A_a(0))\leq A_a^l(0)<A_{a+1}(0)$.
This yields that $A^l_a(0)$ is in $k$ normal form.\hfill$\Box$

\begin{defi} We define the base change operations recursively as follows.
\begin{enumerate}
\item If $m=0$ then $m[k\leftarrow k+1]:=0$.
If $m\gknf A_a(k,b)+l>0$ then 
$m[k\leftarrow k+1]:=A_{a[k\leftarrow k+1]}(k+1,b[k\leftarrow k+1])+l$.
\item  If $m=0$ then $m\scbr{k\leftarrow k+1}:=0$.
If $m\gknf A_a(k,b)+l>0$ then 
$m\scbr{k\leftarrow k+1}:=A_a(k+1,b\scbr{k\leftarrow k+1})+l$.
\item If $m=0$ then $m\{k\leftarrow k+1\}:=0$. If $m\gknf A_a(k,b)+l>0$ then 
$m\{k\leftarrow k+1\}:=A_{a\{k\leftarrow k+1\}}(k+1,b)+l$.
\end{enumerate}
\end{defi}

\begin{lem} Let $m':=m[k\leftarrow k+1]$, $m'':=m\scbr{k\leftarrow k+1}$, and $m''':=m\{k\leftarrow k+1\}.$
\begin{enumerate}
\item $m\leq m'$, $m\leq m''$ and $m\leq m'''$. Moreover
$A_a(k,b)\leq A_{a[k\leftarrow k+1]}(k+1,b[k\leftarrow k+1])$, 
$A_a(k,b)\leq A_{a\scbr{k\leftarrow k+1}}(k+1,b\scbr{k\leftarrow k+1})$,
$A_a(k,b)\leq A_{a\{k\leftarrow k+1\}}(k+1,b\{k\leftarrow k+1\})$ even if $A_a(k,b)$ is not in $k$-normal form.
\item If $m\geq k$ then $m<m'$, $m< m''$ and $m<m'''$.
\item If $m>0$ then $(m-1)'<m'$, $(m-1)''<m''$, and $(m-1)'''<m'''$.
\item If $m\gknf A_a(k,b)+l$ then $m'=_{k+1-\mathrm{nf}}A_{a'}(k+1,b')+l$, $m''=_{k+1-\mathrm{nf}}A_{a}(k+1,b'')+l$,
and $m'''=_{k+1-\mathrm{nf}}A_{a'''}(k+1,b)+l$.
\end{enumerate}
\end{lem}
The first claim of the first assertion is proved by induction on $m$.
Assume that $m\gknf A_a(b)+l$.

Case 1. $a=0$. 

Then $0<m<k$. Then $m=m'$.

Case 2. $a>0$. 
Then  the induction hypothesis yields $m= A_a(b)+l\leq B_{a'}(b')+l=m'$.
The other claims in the first assertion are proved similarly.

The first claim of the second assertion can be seen by inspection of the proof of the first assertion since $A_a(b)<B_a(b)$ for $a>0$.

The first claim of the third assertion is proved by induction on $m$.

Assume that $m\gknf A_a(b)+l$.

Case 1. $a=0$. Then $0<m<k$. Then $(m-1)'=m-1<m=m'$.

Case 2. $a>0$. 

Case 2.1. $l>0$. 

Then $m-1\gknf A_a(b)+l-1$
and  $(m-1)'=A_a(b)'+l-1<A_a(b)'+l=m'$.

Case 2.2. $l=0$.

Case 2.2.1. $b>0$. 

Then for some $p$ we find $m-1\gknf A_a(b-1)+p<m\gknf A_a(b).$

If $a=1$ then $p<k$ hence $p<k+1$ and then the induction hypothesis yields $(m-1)'=A_a(b-1)'+p\leq B_{a'}((b-1)')+p \leq B_{a'}(b'-1)+p<B_{a'}(b')= m'$.

If $a\geq 2$ then    $(m-1)'=A_a(k,b-1)'+p$ where $p<A_a(b)=A_{a-1}^k(A_a(b-1)) \leq B_{a'-1}^k(B_{a'}(b'-1))$ and we arrive at 
$(m-1)'=A_a(k,b-1)'+p<B_{a'-1}^{k}(B_{a'}(b'-1))\cdot 2< B_{a'-1}^{k+1}(B_{a'}(b'-1))=B_{a'}(b')=m'$.

Case 2.2.2. $b=0$. 

Then $m-1=A_a(0)-1=A_{a-1}^{k}(0)-1=A_{a-1}(A_{a-1}^{k-1}(0)-1)+p$ for some $p<\om$.

If  $a=1$ then $(m-1)'=(A_1(0)-1)'=k-1<k+1=A_1(0)'.$  

If $a>1$ then the induction hypothesis yields $(m-1)'=(A_{a-1}^k(0)-1)'
=(A_{a-1}(A_{a-1}^{k-1}(0)-1))'+p$.
We have $p<A_a(0)=A_{a-1}^k(0)\leq B_{a'-1}^k(0)$ 
and the induction hypothesis yields $A_{a-1}(A_{a-1}^{k-1}(0)-1)'\leq B_{(a-1)'}((A_{a-1}^{k-1}(0))')
\leq B_{(a-1)'}^k(0) \leq B_{a'-1}^k(0)$.
Hence $(m-1)'=(A_{a-1}^k(0)-1)'
=A_{a-1}(A_{a-1}^{k-1}(0)-1)'+p<B_{a'-1}^k(0)\cdot 2<B_{a'-1}^{k+1}(0)=B_{a'}(0)$.

The second and third claim in the third assertion are proved by a similar induction on $m$.

Let us prove the first claim in the fourth assertion.
Assume that $A_a(0)\leq m<A_{a+1}(0)$ and 
$A_a(b)\leq m<A_{a}(b+1)$ and that $m=A_a(b)+l$.
Then $b<A_a^{k-1}(0)$ hence $b'<B_{a'}^{k-1}(0)$ by the third assertion.
We have $l<A_a^k(0)\leq B_{a'}^k(0)$. 
Hence $B_{a'}(0)\leq m' \leq B_{a'}^k(0)\cdot 2< B_{a'}^{k+1}(0)=B_{a'+1}(0).$
(Please note that $B_{a'+1}(0)$ might be smaller than $B_{(a+1)'}(0)$.)
$A_a(b)\leq m$ yields $B_{a'}(b')\leq m'=B_{a'}(b')+l$ by the third assertion.
We have $l<A_{a-1}^k(A_a(b))\leq B_{a'-1}^k(B_{a'}(b'))$.
Hence $m'=B_{a'}(b')+l<B_{a'-1}^{k+1}(B_{a'}(b'))=B_{a'}(b'+1)$.
This shows that $m'=B_{a'}(b')+l$ is in $k+1$ normal form.

The second and third claim in the fourth assertion are proved by a similar induction on $m$.
\hfill$\Box$

We define for technical reasons $-m+\al$ in the natural way. 
Thus  $-m+\al$ is equal to $\al$ if $\alpha$ is infinite and we set $-m+\al$ be equal to the maximum
of $\{-m+\al,0\}$ if $\al<\om$. 

\begin{defi}
\begin{enumerate}
\item $\pk0:=0$,
\item If $m\gknf A_a(k,b)+l$ then $$\pk m:=
\left\{
\begin{array}{lr}
b+1& \mbox{if }a=0\\
\om\cdot (1+\pk b)+l&\mbox{if } a=1\\
\om^{\om+(-2+\pk a)}+\om\cdot \pk b+l &\mbox{if } a\geq 2
\end{array}\right\}.$$
\end{enumerate}
\end{defi}
\begin{defi}
\begin{enumerate}
\item $\ck0:=0$
\item  If $m\gknf A_a(k,b)+l$ then  $$\ck m:=
\left\{
\begin{array}{lr}
b+1& \mbox{if }a=0\\
\om\cdot (1+\ck b)+l&\mbox{if } a=1\\
\om^{\om+(-2+ a)}+\om\cdot \ck b+l &\mbox{if } a\geq 2
\end{array}\right\}.$$
\end{enumerate}
\end{defi}
\begin{defi}
\begin{enumerate}
\item $\xk0:=0$
\item  If $m\gknf A_a(k,b)+l$ then  $$\xk m:=
\left\{
\begin{array}{lr}
b+1& \mbox{if }a=0\\
\om\cdot (1+\ck b)+l&\mbox{if } a=1\\
\om^2\cdot (-1+\xk a) +\om\cdot b+l&\mbox{if } a\geq 2
\end{array}\right\}.$$
\end{enumerate}
\end{defi}

Recall that we write $m':=m[k\leftarrow k+1]$, $m'':=m\scbr{k\leftarrow k+1}$, and $m''':=m\{k\leftarrow k+1\}$
when $k$ is fixed in the context.
\begin{lem}
 
\begin{enumerate}
\item $\pskpe m'=\pk m$, $\cskpe m''=\ck m$, and $\xskpe m''=\xk m$.
\item If $m>0$ then $\pk (m-1)<\pk m$, $\ck (m-1)<\ck m$ and $\xk (m-1)<\xk m$.
\end{enumerate}
\end{lem}
The first claim in the first assertion is proved by an easy induction on $m$.
The claim is clear for $m=0$.
Assume that $m\gknf A_a(b)+l$. If $a=0$ then $m<k$ and $\pskpe m'=m=\pk m$.
If $a=1$ then $\pskpe m'=\om\cdot (1+\pskpe b')+l=\om\cdot (1+\pk b)+l=\pk m$.
If $a>1$ then $\pskpe m'=\om^{\om+(-2+\pskpe a')}+\om\cdot \pskpe b'+l=\om^{\om+(-2+\pk a)}+\om\cdot \pk b+l=\pk m$.
The second and third claim in the first assertion follow similarly.

The first claim in the second assertion is proved by induction on $m$.

Assume that $m\gknf A_a(b)+l$.

Case 1. $a=0$. Then $m<k$. Then $\pk (m-1)=m-1<m=\pk m$.

Case 2. $a>0$. 

Case 2.1. $l>0$. $m-1\gknf A_a(b)+l-1$.
Then we find $\pk (m-1)=\pk (A_a(b))+l-1<\pk (A_a(b))+l=\pk m$.

Case 2.2. $l=0$.

Case 2.2.1. $b>0$. Then for some $p$ we find $m-1\gknf A_a(b-1)+p<m\gknf A_a(b).$
If $a=1$ then the induction hypothesis yields $\pk(m-1)=\om\cdot (1+\pk(b-1))+p<\om\cdot (1+\pk b)=\pk m$.

If $a\geq 2$ then the induction hypothesis yields  $\pk(m-1)=\pk (A_a(b))+ \om\cdot \pk(b-1)+p<\pk (A_a(b)) +\om\cdot \pk (b)=\pk m$.

Case 2.2.2. $b=0$.
If  $a=1$ then $\pk(m-1)=\pk(A_1(0)-1)=\pk(k-1)=k-1<\om=\pk(A_1(0))=\pk m.$

If $a=2$ then  for certain $p_1,\ldots,p_k$ we find $\pk(m-1)=\pk(A_1^k(0)-1)=\om\cdot(\pk (A_1^{k-1}(0)-1)+p_1 )=\om^k\cdot p_k+\ldots +\om^1\cdot p_1<\om^\om=\pk m.$

If $a>2$ then for certain $p_0,\ldots,p_{k-1}$ the induction hypothesis yields 

\begin{eqnarray*}
&&\pk(m-1)=\pk(A_{a-1}^k(0)-1)\\
&=&
\pk( A_{a-1}((A_{a-1}^{k-1}(0)-1))+p_0 )\\
&=&\om^{\om+(-2+\pk (a-1))}+\om\cdot  \pk (A_{a-1}^{k-1}(0)-1)+p_0\\
&\leq &\om^{\om+(-2+\pk (a-1))}\cdot k+\om^{k-1}\cdot p_{k-1}+\cdots+\om^0\cdot p_0\\
&<&\om^{\om+(-2+\pk (a))}\\
&=&\pk m.
\end{eqnarray*}\

The second claim in the second assertion is proved by a similar induction on $m$.

The third claim in the second assertion is proved by induction on $m$.
Let us just consider the case $m\gknf A_a(b)+l$ where $a\geq 2$, $b=0$, $l=0$.
Then the induction hypothesis yields 
\begin{eqnarray*}
&&\xk(m-1)=\xk(A_{a-1}^k(0)-1)\\
&=&
\xk( A_{a-1}((A_{a-1}^{k-1}(0)-1))+p_0 )\\
&=&\om^2\cdot (-1+\xk (a-1))+\om\cdot  (A_{a-1}^{k-1}(0)-1)+p_0\\
&<&\om^2(-1+\xk (a))\\
&=&\xk m.
\end{eqnarray*}\

\hfill$\Box$

\begin{defi}
Let $m<\om$. 
\begin{enumerate}
\item
Put $m_0:=m.$
Assume recursively that $m_{l}$ is defined and $m_{l}>0$.
Then $m_{l+1}=m_l[l+3\leftarrow l+4]-1$. If $m_l=0$ then $m_{l+1}:=0$.
\item Put $\tilde{m}_0:=m.$
Assume recursively that $\tilde{m}_{l}$ is defined and $\tilde{m}_{l}>0$.
Then $\tilde{m}_{l+1}=\tilde{m}_l\scbr{l+3\leftarrow l+4}-1$. If $\tilde{m}_l=0$ then $\tilde{m}_{l+1}:=0$.
\item Put $\overline{m}_0:=m.$
Assume recursively that $\overline{m}_{l}$ is defined and $\overline{m}_{l}>0$.
Then $\overline{m}_{l+1}=\overline{m}_l\{l+3\leftarrow l+4\}-1$. If $\overline{m}_l=0$ then $\overline{m}_{l+1}:=0$.
\end{enumerate}
\end{defi}

\begin{theo} 
\begin{enumerate}
\item For all $m<\om$ there exists an $l<\om$ such that $m_l=0.$ This is provable in $\mathrm{PRA}+\mathrm{TI(}\varepsilon_0)$.
\item For all $m<\om$ there exists an $l<\om$ such that $\tilde{m}_l=0.$ This is provable in $\mathrm{PRA}+\mathrm{TI}(\om^{\om+\om})$.
\item For all $m<\om$ there exists an $l<\om$ such that $\overline{m}_l=0.$ This is provable in $\mathrm{PRA}+\mathrm{TI}(\om^{\om})$.
\end{enumerate}
\end{theo} 
Proof. Define $o(m,l):=\psi_{l+3}( m_l).$ If $m_{l+1}>0$ then by the previous lemmata
\begin{eqnarray*}
o(m,l+1)&=&\psi_{l+4}(m_{l+1})   \\
&=&\psi_{l+4} (m_l[l+3\leftarrow l+4]-1)\\
&<&\psi_{l+4} (m_l[l+3\leftarrow l+4])\\
&=&\psi_{l+3} (m_l)\\
&=&o(m,l)
\end{eqnarray*}
This proves the first assertion. The second and assertion are proved similarly by now using $\ck$ ($\xk$ resp.) instead of $\pk$.

Let us now prove the independence results. Recall that the standard system of fundamental sequences for the ordinals less than $\eo$
is defined recursively as follows.
If $\al=0$ then $\al[x]:=0$.
If $\al=\be+1$ then $\al[x]:=\be$.
If $\al=\om^{\al_1}+\cdots +\om^{\al_n}$ where $\al_1\geq \ldots\geq \al_n$ and if $\al_n$ is limit then 
$\al[x]=\om^{\al_1}+\cdots+ \om^{\al_n[x]}$.
If $\al=\om^{\al_1}+\cdots +\om^{\al_n+1}$ where $\al_1\geq \ldots\geq \al_n+1$ then 
$\al[x]=\om^{\al_1}+\cdots +\om^{\al_n}\cdot x$.
Then for limit ordinals $\la<\eo$ we have $\la[x]<\la[x+1]<\la$ and $\la=\sup\{\al[x]:x<\om\}$.
Moreover note that $-2+x=\om[x-2]$ and 
$-2+\la[x]\leq (-2+\la)[x]$ holds for $\la>\om$.

It is easy to show that these fundamental sequences fulfill the so called Bachmann property:
If $\al[x]<\be<\al$ then $\al[x]\leq\be[1]$. (See, e.g., \cite{BCW}, for a proof and further discussion.)

Let $\al\geq_l\be$ iff there exist $\al_0,\ldots,\al_m$ such that $\al=\al_0$, $\be=\al_m$ and $\al_{i+1}=\al_i[l]$ for all $i<m$.
If $\al[x]<\be<\al$ then the Bachmann property yields $\be\geq_1 \al[x]$. Moreover  the Bachmann property yields that $\al\geq_l\be$ implies $\al\geq_{l+1}\be$.

\begin{lem}\label{bach}
Assume that $m>0$.
\begin{enumerate} 
\item $\pk m>\pskpe (m'-1)\geq (\pk m)[k-2]$.
\item $\ck m>\cskpe (m''-1)\geq (\ck m)[k-2]$.
\item $\xk m>\xskpe (m'''-1)\geq (\xk m)[k-1]$.
\end{enumerate}
\end{lem}
Proof. 
Let us first prove the first assertion. Clearly $\pk m=\pskpe (m')>\pskpe (m'-1).$ So let us prove the second inequality.

If $0<m<k$ then $\pskpe(m'-1)=m-1= (\pk m)[k-2].$

Assume that $m\gknf A_a(k,b)+l\geq k$. Then $a>0$.

Case 1. $l>0$. Then $m-1\gknf A_a(k,b)+l-1$.
Then $(\pk m)[k-2]=(\pk (A_a(k,b))+l)[k-2]=\pk (A_a(k,b))+l-1$ and 
$\pskpe (m'-1)=\pskpe(B_{a'}(b')+l-1)=\pk (A_a(k,b))+l-1$.

Case 2. $l=0$.

Case 2.1. $b>0$. 

If $a=1$ then the induction hypothesis yields $(\pk m)[k-2]=(\om(1+\pk b))[k-2]\leq \om(1+(\pk b)[k-2])+k\leq \om(1+(\pskpe (b'-1)))+k= \pskpe(B_{a'}(b'-1)+k\leq \pskpe(B_{a'}(b')-1)=\pskpe (m'-1)$.

If $a>1$ then the induction hypothesis yields
$(\pk m)[k-2]=(\om^{\om+(-2+\pk a)} +\om\cdot \pk b)[k-2] \leq   \om^{\om+(-2+\pk a)} +\om\cdot ((\pk b)[k-2] )+k\leq 
\om^{\om+(-2+\pskpe a')} +\om\cdot (\pskpe (b'-1))+k=\pskpe(B_{a'}(b'-1)+k\leq \pskpe(B_{a'}(b')-1)=\pskpe (m'-1)$.

Case 2.2. $b=0$.
If $a=1$ then $m=A_1(0)=k$ and $(\pk m)[k-2]=\om[k-2]\leq k= \pskpe(k+1-1)=(\pskpe (m'-1)$.

If $a=2$ then $(\pk m)[k-2]=(\om^{\om+(-2+\pk(2)})+\om\cdot 0)[k-2]=\om^\om[k-2]\leq \om^k=\om(1+\pskpe(B_1^{k-1}(0))=\pskpe(B_1^k(0))\leq 
\pskpe (B_1^{k+1}(0)-1)=(\pskpe (m'-1)$. Note that $\pskpe B_1^l(0)=\om^l$ holds by induction on $l$ for $0<l\leq k$.

If $a=k$ then

\begin{eqnarray*}
&&(\pk m)[k-2]=(\om^{\om+(-2+\pk a)})[k-2]\\
&\leq &(\om^{\om+(-2+\om)})[k-2]\\
&= &(\om^{\om+(-2+k)})\\
&=& (\om^{\om+(-2+\pskpe (a'-1))})\cdot k\\
&=&\pskpe(B_{a'-1}^k(0))\\
&\leq &\pskpe(B_{a'}(0)-1)\\
&=&\pskpe (m'-1).
\end{eqnarray*}

If $a>2$ and $a\not=k$ then $(-2+\pk a)[k-2]=-2+(\pk a)[k-2]$ and

\begin{eqnarray*}
&&(\pk m)[k-2]=(\om^{\om+(-2+\pk a)})[k-2]\\
&\leq& (\om^{\om+(-2+(\pk a)[k-2])})\cdot k\\
&\leq& (\om^{\om+(-2+\pskpe (a'-1))})\cdot k\\
&=&\pskpe(B_{a'-1}^k(0))\\
&\leq &\pskpe(B_{a'}(0)-1)\\
&=&\pskpe (m'-1).
\end{eqnarray*}

Note that $\pskpe (B_{a'-1}^l(0))=(\om^{\om+(-2+\pskpe (a'-1)})\cdot l$ holds by induction on $l$ for $0<l\leq k$.

The second assertion is proved by a similar induction on $m$. An even stronger result is proved later, see assertion 2 of Lemma \ref{lem15}.

Let us finally prove the third assertion. Clearly $\xk m=\xskpe (m''')>\xskpe (m'''-1).$ So let us prove the second inequality.

If $0<m<k$ then $\xskpe(m'''-1)=m-1= (\xk m)[k-1].$

Assume that $m\gknf A_a(k,b)+l\geq k$. Then $a>0$.

Case 1. $l>0$. Then $m-1\gknf A_a(k,b)+l-1$.
Then $(\xk m)[k-1]=(\xk (A_a(k,b))+l)[k-1]=\xk (A_a(k,b))+l-1$ and 
$\xskpe (m'''-1)=\xskpe(B_{a}(b)''')+l-1=\xk (A_a(k,b))+l-1$.

Case 2. $l=0$.

Case 2.1. $b>0$. 

If $a=1$ then the induction hypothesis yields $(\xk m)[k-1]=\om(1+ b)[k-1]\leq \om(1+b[k-1])+k\leq \om(1+(b-1))+k= \xskpe(B_{a'''}(b-1)+k\leq \xskpe(B_{a'''}(b)-1)=\xskpe (m'''-1)$.

If $a>1$ then the induction hypothesis yields
$(\xk m)[k-1]=(\om^2(-1+\xk a) +\om\cdot b)[k-1] \leq   (\om^2(-1+\xk a) +\om\cdot b)[k-1] )+k\leq 
\om^2(-1+\xskpe a''') +\om\cdot (b-1))+k=\xskpe(B_{a'''}(b-1)+k\leq \xskpe(B_{a'''}(b)-1)=\xskpe (m'''-1)$.

Case 2.2. $b=0$.
If $a=1$ then $m=A_1(0)=k$ and $(\xk m)[k-1]=\om[k-1]\leq k=  \xskpe(k+1-1)=(\xskpe (m'''-1)$.

If $a=2$ then $(\xk m)[k-1]=(\om^2(-1+\xk(2))+\om\cdot 0)[k-1]=\om^2[k-1]\leq \om\cdot k\leq \om(1+(B_1^{k-1}(0))=\xskpe(B_1^k(0))\leq 
\xskpe (B_1^{k+1}(0)-1)=(\xskpe (m'''-1)$. Note that $\xskpe (B_1^l(0))\geq \om\cdot l$ holds by induction on $l$ for $0<l\leq k$.

If $a=k$ then

\begin{eqnarray*}
&&(\xk m)[k]=(\om^2(-1+\xk a))[k-1]\\
&\leq &(\om^3)[k-1]\\
&= &\om^2\cdot (k-1)\\
&=& \om^{2}(-1+\xskpe (a'''-1))\\
&\leq &\xskpe(B_{a'''-1}(0))\\
&\leq &\xskpe(B_{a'''}(0)-1)\\
&=&\xskpe (m'''-1)
\end{eqnarray*}

If $a>2$ and $a\not=k$ then 

\begin{eqnarray*}
&&(\xk m)[k-1]=\om^2(-1+\pk a))[k-1]\\
&\leq& \om^2(-1+(\xk a)[k-1]))+\om \cdot k\\
&\leq& \om^2(-1+\pskpe (a'''-1))+\om\cdot B_{a'''-1}^{k-1}(0)\\
&=&\xskpe(B_{a'''-1}^k(0))\\
&\leq &\xskpe(B_{a'''}(0)-1)\\
&=&\xskpe (m'''-1)
\end{eqnarray*}

\hfill$\Box$

For $\al<\eo$ let $\mathrm{PRA}+\mathrm{TI}(<\!\al)$ be the union of the theories $\mathrm{PRA}+\mathrm{TI}(\be)$ where $\be<\al$.

\begin{theo}
\begin{enumerate}
\item $\mathrm{PA}\not \vdash (\forall m)( \exists l )[m_l=0]$.
\item $\mathrm{PRA}+\mathrm{TI}(<\om^{\om+\om})\not \vdash (\forall m )(\exists l)[ \tilde{m}_l=0]$.
\item $\mathrm{PRA}\not \vdash (\forall m)( \exists l)[ \overline{m}_l=0]$.
\end{enumerate}
\end{theo}

Proof of the first assertion. Let $m(1):=A_2(3,0)$ and $m(r+1):=A_{m(r)}(3,0)$. Let $\om_1:=\om$ and $\om_{r+1}:=\om^{\om_r}$.
Then $\psi_3(m(r))=\om_{r+1}$ for $r\geq 1$.
We claim that $o(m(r),l)\geq_1\om_{r+1}[1]\ldots [l]$.
Proof of the claim. Write $m$ for $m(r)$.
For $o(m,l)>0$ Lemma \ref{bach} yields $o(m,l)>o(m,l+1)=\psi_{l+4}(m_l[l+3\leftarrow l+4]-1)\geq (\psi_{l+3}(m_l))[l+1]=(o(m,l)[l+1]$.
The Bachmann property yields $o(m,l+1)\geq_1 (o(m,l)[l+1]$.
The induction hypothesis yields $o(m,l)\geq_1 \om_{r+1}[1]\ldots [l]$.
This yields $o(m,l)[l+1]\geq_1 (o(m,l)[l+1]\geq_1 \om_{r+1}[1]]\ldots [l][l+1]$.

Therefore the least $l$ such that $o(m,l)=0$ is at least as big as the least $l$ such that 
$\om_{r+1}[1]\ldots [l+1]=0$. The result follows from $\mathrm{PA}\not\vdash \forall r \exists l (\om_{r+1})[2]\ldots [l+1]=0$.

The second assertion follows similarly. 

Let $m((r)):=A_{r+2}(3,0)$. Then $\chi_3(m(r))=\om^{\om+r}$.
The result follows from $\mathrm{PRA} \not\vdash \forall r \exists l (\om^{\om+r})[1]\ldots [l]=0$.

The third assertion follows similarly by employing $\xk$ and $m(r)$.\hfill$\Box$

It is somewhat surprising that one obtains a $\mathrm{PA}$ independence via $\forall m \exists l m_l=0$
without using the nested sandwiching procedure from \cite{Weiermann}. We consider this Goodstein principle for the Ackermann function
as a giant Goodstein principle.

Note that the assertion $\forall m \exists l \tilde{m}_l=0$ is considerably weaker than the corresponding
assertion in \cite{Weiermann} where the base change in the second argument of the Ackermann function  led
to independence from $\mathrm{I}\Sigma_2$. We consider this Goodstein principle therefore as an intermediate Goodstein principle.

The assertion $\forall m \exists l \overline{m}_l=0$ is an example for an illusionary giant Goodstein principle.
In strength it does not exceed the axioms needed for proving the totality of the functions which are involved in the definition
of the underlying normal forms. But this principle is still non trivial since it is independent of $\mathrm{PRA}$.

An even weaker Goodstein principle can be obtained by performing a trivial base change in the iteration parameter $k$.
For $m\gknf A_a(k,b)+l\geq k$ let $m'''':=A_a(k+1,b)+l$. A Goodstein principle base on this definition becomes even provable in
in a weak meta theory using the ordinal assignment $ord(m):=\om^2\cdot a+\om b+l$. Note that the graph of the Ackermann function is
elementary and so it makes sense to speak about writing $m$ in $k$ normal form even in a very weak base theory. But defining the base change
already requires some tricky machinery.

\section{Giant and illusionary giant Goodstein sequences for the extended Grzegorczyk hierarchy}
Let us agree that small Greek letters denote ordinals less than $\eo$.

Let us recall the definition of the standard assignment of fundamentals sequences for the ordinals below $\varepsilon_0$.
We agree on $\om^{\al+1}(\be+1)[x]=\om^{\al+1}\be+\om^\al\cdot x$
and $\om^{\la}(\be+1)[x]=\om^{\la}\be+\om^{\la[x]}.$
We also assume for convenience that $(\al+1)[x]:=\al$ and $0[x]:=0$.

Using these fundamental sequences we define the extended Ackermann functions $\al,b\mapsto A_\al(k,b)$ for $\al<\eo$ and $k,b<\om$ recursively as follows:

\begin{eqnarray*}
A_0(k,b)&:=&b+1,\\
A_{\al+1}(k,0)&:=&A_\al(k,\cdot)^k(0),\\
A_{\al+1}(k,b+1)&:=&A_\al(k,\cdot)^k(A_{\al+1}(k,b)),\\
A_{\la}(k,0)&:=&A_{\la_{k,k,0}}(k,\cdot)^k(0),\\
A_{\la}(k,b+1)&:=&A_{\la_{k,k,A_{\la}(k,b)}}(k,\cdot)^k(A_{\la}(k,b)).
\end{eqnarray*}
where for a limit $\la$ the ordinal $\la_{l,k,b}$ is defined recursively by $\la_{0,k,b}:=\la[b]$ and $\la_{l+1,k,b}:=\la[A_{\la_{l,k,b}}(k,b)]$.

For the rest of this article by $k$ will denote a positive integer not smaller than $3$.

In the sequel we often write $A_\al(b)$ for $A_\al(k,b)$ and $B_\al(b)$ for $A_\al(k+1,b)$ and
$\la_{l,b}$ for $\la_{l,k,b}$ when $k$ is fixed in a given context.

The functions $A_\al$ come along with natural monotonicity properties due to the Bachmann property of the system of fundamental sequences.
Recall that his property states that $\al[l]<\be<\al$ yields $\al[l]\leq \al[1]$. 
Moreover recall that $\leq_l$ is the transitive and reflexive closure of $\{(\al[l],\al):\al<\eo\}$. 

In general the function $\al\mapsto A_\al(b)$ is not monotone in $\al$ but it shows decent monotonicity with respect to the relation $\leq_l$.

Furthermore, for $\al=\om^{\al_1}\cdot m_1+\cdots+\om^{\al_n}\cdot m_n$ with $\al_1>\ldots> \al_n$ and $0<m_1,\ldots,m_n<\om$ let 
$mc(\al):=\max\{mc(\al_1),\ldots,mc(\al_n),m_1,\ldots,m_n\}$ where $mc(0):=0$. We call $mc(\al)$ the maximal coefficient  of an ordinal $\al$.
This maximal coefficient plays an important role in bounding values of $A_\al(b)$.

\begin{lem}\label{standard}
\begin{enumerate}
\item  $A_\al(b)<A_\al(b+1)$,
\item $\al[l]<\be<\al$ yields $A_{\al[l]}(b)<A_\be(b)$,
\item $\al\leq_l\be$ yields $A_\al(b)\leq A_\be(b)$ for all $b\geq l\geq 1$,
\item $mc(\al)<A_\al(b)$,
\item $\al<\be$ and $mc(\al)\leq b$ yields $A_\al(b)\leq A_\be(b)$.
\end{enumerate}
\end{lem}
Proof. This is easy. One can e.g. consult \cite{BCW} or \cite{Weiermannb} if needed.\hfill$\Box$

\begin{lem}
For any $m$ there will be at most finitely many $\al<\eo$ such that $A_\al(k,0)\leq m$.
\end{lem}
Proof. First note that for a given $m$ there will be finitely many $\al<\eo$ such that $N\al\leq m$
where $N\al$ is the number of occurrences of $\om$ in the Cantor normal form of $\al$ (cf., e.g., \cite{BCW}).
An easy induction on $\al$ yields that $A_\al(k,b)\geq N\al+b$.
Putting things together the Lemma follows.\hfill $\Box$

\begin{lem}
For all $m>0$ there exist unique $\al<\varepsilon_0$ and $b,l<\om$ such that
\begin{enumerate}
\item $m=A_\al(k,b)+l,$
\item $\al$ is maximal with $A_\al(k,0)\leq m$ (so that $A_{\al}(k,0)\leq m <A_{\al+1}(k,0)$),
\item $b$ is maximal with $A_\al(k,b)\leq m$ (so that $A_\al(k,b)\leq m <A_{\al}(k,b+1)).$
\end{enumerate}
\end{lem}
We write $m\gknf A_\al(k,b)+l$ in this case and call $A_\al(k,b)+l$ the $k$ normal form of $m$.
This normal form is uniquely determined.
If $m\gknf A_\al(k,b)+l$ and $a=0$ then necessarily $m=A_0(k,\cdot)^m(0)$ where $m<k$. 
If $m=k$ then $m\gknf A_1(k,0)$ and $A_1(k,b)=k\cdot (1+b)$.
The yields in particular that $A_\al(k,b)> 2\cdot b$ for $\al>0$.
The latter estimate is used tacitly at several occasions.

\begin{lem}
\begin{enumerate}
\item If $m\gknf A_\al(b)$ and $b>0$ then $A_\al(b-1)$ is in $k$ normal form, too.
\item $A^l_\al(0)$ is in $k$ normal form for $0<l<k$.
\item If $m\gknf A_\la(0)$ and $\la$ is a limit then $A_{\la_{l,0}}(0)$ is in $k$ normal form for $0<l<k$.
\item If $m\gknf A_\la(b)$ and $\la$ is a limit and $b>0$ then $A_{\la_{l,b}}(A_\la(b-1))$ is in $k$ normal form for $0<l<k$.
\end{enumerate}
\end{lem}
Proof. This follows easily from the Bachmann property.\hfill$\Box$

\begin{defi} We define the base change operations recursively as follows.
\begin{enumerate}
\item If $m=0$ then $m[k\leftarrow k+1]:=0$.
If $m\gknf A_\al(k,b)+l>0$ then 
$m[k\leftarrow k+1]:=A_{\al[k\leftarrow k+1]}(k+1,b[k\leftarrow k+1])+l$.
If $\al=\om^\be\cdot m+\ga$ is in Cantor normal form then
$\al[k\leftarrow k+1]=\om^{\be[k\leftarrow k+1]}\cdot m[k\leftarrow k+1]+\ga[k\leftarrow k+1].$ 

\item  If $m=0$ then $m\scbr{k\leftarrow k+1}:=0$.
If $m\gknf A_\al(k,b)+l$ then 
$m\scbr{k\leftarrow k+1}:=A_\al(k+1,b\scbr{k\leftarrow k+1})+l$.

\item  If $m=0$ then $m\{k\leftarrow k+1\}:=0$.
If $m\gknf A_\al(k,b)+l$ then 
$m\{k\leftarrow k+1\}:=A_{\al\{k\leftarrow k+1\}}(k+1,b)+l$.
If $\al=\om^\be\cdot m+\ga$ is in Cantor normal form then
$\al\{k\leftarrow k+1\}=\om^{\be\{k\leftarrow k+1\}}\cdot m\{k\leftarrow k+1\}+\ga\{k\leftarrow k+1\}.$ 
\end{enumerate}
\end{defi}

It is highly non trivial to show that the base change operations preserve monotonicity and 
normal forms. To show these properties we develop some new machinery about fundamental sequences.

For an ordinal context $\la[[\cdot ]]$ with exactly one occurrence of the placeholder $[[\cdot ]]$ we define an ordinal context $\la^*[[\cdot ]]$, the truncation of  $\la[[\cdot ]]$ as follows.
If  $\la[[\cdot ]]=[[\cdot ]]$ then $\la^*[[\cdot ]]=[[\cdot ]]$.
If $\la[[\cdot ]]=\om^{\al_1}+\cdots+ \om^{\al_{i}}+[[\cdot ]]+\om^{\al_{i+1}}+\cdots+\om^{\al_n}$ then 
$\la^*[[\cdot ]]=\om^{\al_1}+\cdots+ \om^{\al_{i}}+[[\cdot ]]$.
If $\la[[\cdot ]]=\om^{\al_1}+\cdots+ \om^{\al_{i}[[\cdot]]}+\om^{\al_{i+1}}+\cdots+\om^{\al_n}$ then
$\la^*[[\cdot ]]=\om^{\al_1}+\cdots+ \om^{\al^*_{i}[[\cdot ]]}$.
So we basically cut off hereditarily terms after the placeholder. (We tacitly assume here as usual that  $\al_1\geq \ldots\geq \al_n$.)

\begin{lem}\label{star}
If $\al<\be$ then $\be=\al+1$ or $\al< \be[1]$ or there exists a context $\lambda$, an ordinal $\ga$ and a natural number $r$ such that $\al=\la[[\om^\ga\cdot r]]$ and
$\be=\la^*[[\om^{\ga+1}]]$. Moreover we have $\al<\la^*[[\om^{\ga}\cdot(r+1)]]$.
\end{lem}

Proof. Assume that $\al=\om^{\al_1}+\cdots+\om^{\al_m}$ and that $\be=\om^{\be_1}+\cdots+\om^{\be_n}$ are both written in Cantor normal form.

Case 1. Assume $m<n$ and $\al_i<\be_i$ for all $i\leq m$. Then of course  $\be=\al+1$ or $\al< \be[1]$.

Case 2. There exists an $i\leq \min\{m,n\}$ such that $\al_i<\be_i$ and $\forall j<i (\al_j=\be_j)$.
If $n>i$ then of course $\al< \be[1]$.
So assume $i=n$. Let us write $\al=\om^{\al_1}+\cdots+\om^{\al_i}\cdot s+\xi$ with $\xi<\om^{\al_i}$.

Case 2.1.
Assume that $\be_i=\de+1$. 

If $\de>\al_i$  then of course $\al< \be[1]$.
If $\de=\al_i+1$ then put $\lambda:=\om^{\al_1}+\cdots+\om^{\al_{i-1}}+[[\cdot]]+\om^{\al_{i+1}}+\cdots+\om^{\al_m}$,
$\ga=\al_i$ and $r:=s$. Then $\al=\la[[\om^\ga\cdot r]]$, $\be=\la^*[[\om^{\ga+1}]]$ and $\al<\la^*[[\om^\ga\cdot (r+1)]]$.

Case 2.2. Assume $\be_i\in Lim$. We have $\al_i<\be_i$. The case $\be_i=\al+1$ is impossible and so
by induction hypothesis there are two cases.

Case 2.2.1. $\al_i< \be_i[1]$. Then $\al<\be[1]$.

Case 2.2.2. $ \be_i=\mu^*[[\om^{\ga+1}]]
$ and $ \al_i=\mu[[\om^{\ga}\cdot r ]]
$ where $\al_i<\mu^*[[\om^\ga\cdot (r+1)]].$

Let $\la:=\om^{\al_1}+\cdots+\om^{\al_{i-1}}+\om^{\mu[[\cdot ]]}+\om^{\al_i}\cdot (s-1) +\om^{\al_{i+1}}+\cdots+\om^{\al_m}$.
Then $\al=\la[[\om^\ga\cdot r]]$, $\be=\la^*[[\om^{\ga+1}]]$ and
$\al<\la^*[[\om^\ga\cdot (r+1)]]$.\hfill$\Box$

For notational reasons we agree on $A_\al(k,-1):=0$.
\begin{lem}\label{maj}
Assume  $\al<\be$. Moreover, assume for all $\lambda$, $\ga$ and $r$: if $\al=\la[[\om^\ga\cdot r]]$ and $\de:=\la^*[[\om^{\ga+1}]]$ then $r< A_{\de_{k-1,A_\de(b-1)}}(A_\de(b-1))$.
Then we obtain either $\al+1=\be$ or
$\al+1\leq  \be[A_{\be_{{k-1,A_\be(b-1)}}}(A_\be(b-1))]$. Moreover $A_{\al+1}(b)\leq A_\be(b)$.
\end{lem}

Proof. 
The assertion is clear for $\be=\al+1$. 
If $\al<\be[1]$ then of course $\al+1\leq  \be[A_{\be_{k-1,A_\be(b-1)}}(A_\be(b-1))]$.
Finally assume by Lemma \ref{star} that 
$\al=\la[[\om^\ga\cdot r]]$ and $\be=\la^*[[\om^{\ga+1}]]$.
Then $\be[l]=\la^*[[\om^\ga\cdot l]]$ for every $l$
and $\be_{k,A_\be(b-1)}=\be[A_{(\la^*[[\om^{\ga+1} ]])_{k-1,A_\be(b-1)}}(A_\be(b-1))   ]$.
By assumption we obtain $$r< A_{(\la^*[[\om^{\ga+1} ]])_{k-1,A_\be(b-1)}}(A_\be(b-1)).$$
Lemma \ref{star} yields 

\begin{eqnarray*}
\al&=&\la[[\om^\ga\cdot r]]\\
&<&\la^*[[\om^\ga\cdot (r+1)]]\\
&\leq &\la^*[[\om^\ga\cdot 
A_{(\la^*[[\om^{\ga+1} ]])_{k-1,A_\be(b-1)}}(A_\be(b-1)) ]]\\
&=&\be[A_{(\la^*[[\om^{\ga+1} ]])_{k-1,A_\be(b-1)}}(A_\be(b-1))]\\
&=&\be[A_{\be_{k-1,A_\be(b-1)}}(A_\be(b-1))].
\end{eqnarray*}

The second claim follows from the first claim by induction on $\be$.
Indeed, the claim is obvious when $\be$ is $\al+1$. Otherwise $\al<  \be[A_{\be_{k-1,A_\be(b-1)}}(A_\be(b-1))]$ and
the induction hypothesis yields $A_{\al+1}(k)\leq A_{\be[A_{\be_{k-1,A_\be(b-1)}}(A_\be(b-1))]}(A_\be(b-1))\leq A_\be(b)$.
\hfill$\Box$

To prove the preservation of normal forms after performing a base change operation the following Lemma will be of key importance (together with Lemma \ref{maj}).

\begin{lem}\label{max}
Assume that there is no $\de>\al$  with $A_\de(b)\leq A_\al(b).$
Then for all $\lambda$, $\ga$ and $r$: if $\al=\la[[\om^\ga\cdot r]]$
and $\be=\la^*[[\om^{\ga+1} ]]$ 
 then $r\leq  A_{\be_{k-1,A_\be(b-1)}}(A_\be(b-1))$.
\end{lem}
Proof. Assume that $\al=\la[[\om^\ga\cdot r]]$. Assume for a contradiction that $$r> A_{(\la^*[[\om^{\ga+1} ]])_{k-1,A_\be(b-1)}}(A_\be(b-1)).$$
Let $\be:=\la^*[[\om^{\ga+1}]]$. Then $\be>\al$. But

\begin{eqnarray*}
&&A_\be(b)\\
&=&A^k_{\be_{k,A_\be(b-1)}}(A_\be(b-1))\\
&=&A^k_{\la^*[[  \om^\ga\cdot A_{(\la^*[[\om^{\ga+1} ]])_{k-1,A_\be(b-1)}}]]}(A_\be(b-1))\\
&<& A_{\la^*[[  \om^\ga\cdot r]]}(b)\\
&\leq& A_{\la[[  \om^\ga\cdot r]]}(b)\\
&=&A_\al(b).
\end{eqnarray*}
where we used that $\la^*[[  \om^\ga\cdot r]]\preceq_1 \la[[  \om^\ga\cdot r]]$.
This contradicts the maximality property of $\al$.\hfill$\Box$

\begin{lem} Let $m':=m[k\leftarrow k+1]$, $m'':=m\scbr{k\leftarrow k+1},$ and $m''':=m\{k\leftarrow k+1\}.$
\begin{enumerate}
\item $m\leq m'$, $\al\leq\al'$, $m\leq m''$, $m\leq m'''$ and $\al\leq\al'''$.
\item $A_\al(b)'\leq B_{\al'}(b')$,  $A_\al(b)''\leq B_{\al''}(b'')$ and $A_\al(b)'''\leq B_{\al'''}(b''')$ even if $A_\al(b)$ is not in normal form.
\item If $m\geq k$ then $m< m'$, $m<m''$, $m< m'''$.
\item If $m>0$ then $(m-1)'<m'$, $(m-1)''<m''$, and $(m-1)'''<m'''$
\item If $m\gknf A_\al(k,b)+l$ then $m'=_{k+1-\mathrm{nf}}A_{\al'}(k+1,b')+l$, $m''=_{k+1-\mathrm{nf}}A_{\al}(k+1,b'')+l$,
and $m'''=_{k+1-\mathrm{nf}}A_{\al'''}(k+1,b)+l$.
\end{enumerate}
\end{lem}
The first three assertions are easy to prove using the last two assertions of Lemma \ref{standard}.

The first claim of the fourth assertion is proved by induction on $m$.

Assume that $m\gknf A_\al(b)+l$.

Case 1. $\al=0$. Then $0<m<k$. Then $(m-1)'=m-1<m=m'$.

Case 2. $\al>0$. 

Case 2.1. $l>0$. $m-1\gknf A_\al(b)+l-1$.
Then we find $(m-1)'=A_\al(b)'+l-1<A_\al(b)'+l=m'$.

Case 2.2. $l=0$.

Case 2.2.1. $b>0$. Then we find $m-1\gknf A_\al(b-1)+p<m\gknf A_\al(b).$
Case 2.2.1.1. $\al=1$.
Then $p<k$ hence $p<k+1$ and then the induction hypothesis yields $\pk(m-1)=A_\al(b-1)'+p\leq B_{\al'}((b-1)')+p \leq B_{\al'}(b'-1)+p<B_{\al'}(b'))= m'$.

Case 2.2.1.2. $\al\geq 2$. We find    $(m-1)'=A_\al(k,b-1)'+p$.

By induction on $l$ one verifies $(\al_{l,A_\al(b-1)})'\leq \al'_{l,k+1,B_{\al'}(b-1)}$ using $(\al[x]'\leq \al'[x']$.
Moreover an induction on $l>0$ yields  $mc((\al_{l,A_\al(b-1)})')\leq B_{\al'_{l-1,k+1,B_{\al'}(b'-1))}}(B_{\al'}(b'-1))$. Moreover we have
$mc((\al_{0,A_\al(b-1)})')\leq B_{\al'}(b'-1)$. Then Lemma \ref{standard} yields

 $p<A_\al(b)=A_{\al_{k,A_\al(b-1)}}^k(A_\al(b-1)) \leq B_{\al_{k,A_{\al}(b-1)}'}^k(B_{\al'}(b-1)')\leq B_{\al'_{k+1,k+1,B_\al'(b'-1))}}^k(B_{\al'}(b'-1))
$

and we arrive at $(m-1)'=A_\al(k,b-1)'+l< B_{\al'_{k+1,k+1,B_{\al'}(b'-1) }}^{k+1}(B_{\al'}(b'-1))=B_{\al'}(b')$.

Case 2.2.2. $b=0$.

If  $a=1$ then $(m-1)'=(A_1(0)-1)'=k-1<k+1=A_1(0)'.$  

If $a>1$ then the induction hypothesis yields for some $p$ that $(m-1)'=(A_{\al}^k(0)-1)'
=A_{\al_{k,0}}(A_{\al_{k,0}}^{k-1}(0)-1)'+p$.
We have $p<A_\al(0)=A_{\al_{k,0} }^k(0)\leq B_{(\al_{k,0})'}^k(0)\leq B_{\al'_{k,k+1,0}}^k(0)$ 
and $A_{\al_{k,0}}(A_{\al_{k,0}}^{k-1}(0)-1)'\leq B_{\al'_{k+1,k+1,0}}^k(0)$.
Hence $(m-1)'=(A_{\al_{k,0}}^k(0)-1)'
=(A_{\al_{k,0}}(A_{\al_{k,0}}^{k-1}(0)-1))'+p<B_{\al'_{k+1,k+1,0}}^k(0)\cdot 2<B_{\al'_{k+1,k+1,0}}^{k+1}(0)=B_{\al'}(0)$.

Let us now prove the first claim in the fifth assertion.

Assume that $m\gknf A_{\al}(b)+l$.  Then $A_{\al}(0)$ is in $k$ normal form. Then there is no $\de>\al$ with $A_{\de}(0)\leq A_{\al}(0)$.
We obtain by Lemma \ref{max} that for all contexts $\la$, ordinals $\ga$ and natural numbers $r$:
If $\al=\la[[\om^\ga\cdot r]]$ and $\be=\la^*[[\om^{\ga+1}]]$ then $r\leq A_{\be_{k-1,0}}(0)$.

Assume now that $\al'=\tilde{\la}[[\om^{\tilde{\ga}}\cdot \tilde{r}]]$.
Then there exist $\la,\ga,r$ such that $\tilde{\la}=\la'$, $\tilde{\ga}=\ga'$ and $\tilde{r}=r'$ and $\al=\la[[\om^\ga\cdot r]]$.
Let ${\be}:=\la^*[[\om^{\ga+1}]]$ and $\tilde{\be}=(\la^*)'[[\om^{\ga'+1}]]$.
To apply Lemma \ref{maj}
we have to show $r'< B_{\tilde{\be}_{k,k+1,0}}(0)$.

We obtain $\tilde{\be}[k']=(\la^*)'[[\om^{\ga'}\cdot k']]=(\be[k])'$.
This yields by induction on $l<k$ that $\tilde{\be}_{l,k+1,0} =(\be_{l,0})'$. 
Indeed we find $\tilde{\be}_{0,k+1,0}=(\la^*)'[[\om^{\ga'}\cdot (0)']]=(\be_{0,0})'$.
Moreover $\tilde{\be}_{l+1,k+1,0}=(\la^*)'[[\om^{\ga'}\cdot B_{\tilde{\be}_{l,k+1,{0}(0)}}]]=
(\la^*[[\om^{\ga}\cdot A_{\be_{l,{0}(0)}}]])'=(\be_{l+1,0})'.$
 As desired $r\leq A_{\be_{k-1,0}}(0)$ yields $$r'\leq (A_{\be_{k-1, 0}}(0))'\leq B_{(\be_{k-1, 0})'}(0)< B_{\tilde{\be}_{k,k+1, 0}}(0).$$

We claim that there is no $\de>\al'$ with $B_\de(0')\leq B_{\al'}(b')+l$.
Assume $\al'<\de$. We claim that $B_\de(0)>m'$. 
Indeed we find $m'=B_{\al'}(b')+l\leq B_{\al'}(B_{\al'}^{k-1}(0))+B_{\al'}^k(0)\leq B_{\al'}^k(0)\cdot 2\leq B_{\al'}^{k+1}(0)=
B_{\al'+1}(0)\leq B_\de(0)$ by Lemma \ref{maj}.
So $\al'$ fulfills the maximality condition and $m'$ is in $k+1$-normal form.

\hfill$\Box$

For denoting ordinals below the Howard Bachmann ordinal we use Buchholz's $\psi$ function from \cite{Buchholz}.

Let $\Om$ be the first uncountable ordinal. Let $C(\al)$ be the least set $C$ such that
\begin{enumerate}
\item $\{0,\Om\}\subset C$,
\item If $\be=\be_1+\ldots+\be_n$  and if $\be_1,\ldots,\be_n\in C$ are additive principal and if
$\be_1\geq \ldots\geq \be_n$ then $ \be\in C$,
\item If $\be=\Om^\ga\cdot \de+\eta$ and $\ga,\de,\eta\in C$ then $\be\in C$.
\item $\be\in C\cap \al \Rightarrow \psi\be\in C$.
\end{enumerate}
Let $\psi \al:=\min\{\xi:\xi\not \in C(\al)\}$.
Then $\psi\al<\Om.$
Moreover $\psi\al\in Lim$ for $\al>0$, $\psi0=1$  and $\psi(\al+1)=\psi\al\cdot \om$.
Then $\psi\om=\om^\om$.
We write $\be=_{NF}\psi \al$ if $\be=\psi \al$ and $\al\in C(\al)$.
Then $\be=_{NF}\psi \al$ and  $\de=_{NF}\psi \ga$ and $\al<\ga$ yield $\be<\de$.

Let $OT$ be defined as follows.

\begin{enumerate}
\item $\{0,\Om\}\subset OT$,
\item If $\be=\be_1+\ldots+\be_n$  and if $\be_1,\ldots,\be_n\in OT$ are additive principal and if
$\be_1\geq \ldots\geq \be_n$ then $ \be\in OT$,
\item If $\be=\Om^\ga\cdot \de+\eta$ and $\ga,\de,\eta\in OT$ then $\be\in OT$.
\item If $\be\in OT$ and $\be\in C(\be) $ then $\psi\be\in OT$.
\end{enumerate}

It is well known that $OT\cap \Om=\psi \varepsilon_{\Om+1}$.
Let $\eta_0:=\psi \varepsilon_{\Om+1}$ be an abbreviation for the Howard Bachmann ordinal.

For the termination proof it will be essential to work with $\psi$ terms in normal form.
For $\al\in OT$ define $G\al\subseteq OT$ as follows.

Let $G0:=\emptyset$, $G(\Om^\al\cdot \be+\ga):=G\al\cup G\be\cup G\ga$ and
$G\psi\al:=G\al\cup \{\al\}$. Then $\psi\al$ is in $\psi$ normal form iff $\al\in C(\al)$ iff $G\al<\al$.
Therefore $G$ can be used to single out normal forms and to prove that $OT$ is a primitive recursive set.
By restricting to terms in $\psi$ normal form there will be no chains of terms in $OT$ like $(\psi\Om,\psi\psi\Om,\psi\psi\psi\Om,\ldots)$.

To deal with the third Goodstein principle for the extended Ackermann function we work also with a slight variant $OT'$ of $OT$.

To define it let us first consider a modification of $C(\al)$.
Let $C'(\al)$ be the least set $C$ such that
\begin{enumerate}
\item $\{0,\Om\}\subset C$,
\item  If $\al=\be_1+\ldots+\be_n+\om\cdot b+l$  and if $\be_1,\ldots,\be_n\geq \om^2\in OT$ are additive principal and if
$\be_1\geq \ldots\geq \be_n$ and $b,l<\om$  then $\al\in C$
\item If $\al=\Om^\be\cdot \ga+\de$ and $\be,\ga,\de\in C$ and $\de<\Om^\be$ and $\ga<\Om$ then $\al\in C$.
\item $0<\be\in C\cap \al \Rightarrow \Psi\be\in C$.
\end{enumerate}
For $\al>0$ let $\Psi \al:=\min\{\xi:\xi\not \in C(\al)\}$.
Then $\Psi\al<\Om.$
Moreover $\Psi\al\in Lim$, $\psi 1=\om^2$  and $\Psi(\al+1)=\Psi\al+ \om^2$.
We write $\be=_{NF}\Psi \al$ if $\be=\Psi \al$ and $\al\in C'(\al)$.
Then $\be=_{NF}\Psi \al$ and  $\de=_{NF}\Psi \ga$ and $\al<\ga$ yield $\be<\de$.

\begin{enumerate}
\item $\{0,\Om\}\subset OT'$,
\item If $\al\geq \om^2$ and $b,l<\om$  then $\al+\om\cdot b+l\in OT'$
\item If $\al=\Om^\be\cdot \ga+\de$ and $\be,\ga,\de\in OT'$ and $\de<\Om^\be$ and $\ga<\Om$ then $\al\in OT'$.
\item If $0<\be\in OT'$ and $\be\in C'(\be) $ then $\Psi\be\in OT'$.
\end{enumerate}

Again it is well known $OT'\cap \Om=\Psi \varepsilon_{\Om+1}$ Moreover $\eta_0=\Psi \varepsilon_{\Om+1}$. (See, for example, Buchholz's contributions in \cite{Buchholza} for a proof.)

\begin{defi}
\begin{enumerate}
\item $\pk0:=0$,
\item If $m\gknf A_\al(k,b)+l$ then $$\pk m:=
\left\{
\begin{array}{lr}
b+1& \mbox{if }\al=0\\
\om(1+\pk b)&\mbox{if } \al= 1\\

\psi(\om+(-2+{\pk \al}))+\om \pk b+l &\mbox{if } \al \geq 2
\end{array}\right\}.$$
Moreover for $\al=\om^\be\cdot m+\ga$ in normal form put $\pk\al:=\Om^{\pk\be}\cdot \pk m+\pk \ga$.
\end{enumerate}
\end{defi}
\begin{defi}
\begin{enumerate}
\item $\ck0:=0$
\item  If $m\gknf A_\al(k,b)+l$ then  $$\ck m:=
\left\{
\begin{array}{lr}
b+1& \mbox{if }a=0\\
\om(1+ \ck b)&\mbox{if } \al= 1\\
\om^{\om+(-2+ \al)}+\om\cdot \ck b+l &\mbox{if } \al\geq 2
\end{array}\right\}.$$
\end{enumerate}

\end{defi}
\begin{defi}
\begin{enumerate}
\item $\xk0:=0$
\item  If $m\gknf A_\al(k,b)+l$ then  $$\xk m:=
\left\{
\begin{array}{lr}
b+1& \mbox{if }a=0\\
\om(1+ b)&\mbox{if } \al= 1\\
 \Psi(-1+ \xk\al)+\om\cdot b+l &\mbox{if } \al\geq 2
\end{array}\right\}.$$
\end{enumerate}
Moreover for $\al=\om^\be\cdot m+\ga$ in normal form put $\xk\al:=\Om^{\xk\be}\cdot \xk m+\xk \ga$.
\end{defi}

In the sequel we write $m':=m[k\leftarrow k+1]$, $m'':=m\scbr{k\leftarrow k+1}$ and $m''':=m\{k\leftarrow k+1\}$ 
when $k$ is fixed in the context.

\begin{lem}
 
\begin{enumerate}
\item $\pskpe m'= \pk m$, $\cskpe m''= \ck m$, and $\xskpe m''= \xk m$,
\item If $m>0$ then $\pk (m-1)<\pk m$, $\ck (m-1)<\ck m$, and $\xk (m-1)<\xk m$.
\item 
If $m=A_\al(b)+l$ and $\al>1$ then $\pk m=\psi(\om+(-2+{\pk \al})+\om\cdot \pk b+l$ is in $\psi$ normal form and in Cantor normal form.
Similarly $\xk m=\Psi(-1+{\xk \al})+\om\cdot \pk b+l$ is in $\Psi$ normal form and in Cantor normal form for $\al>1$.
\end{enumerate}
\end{lem}
The first assertion is proved by induction on $m$.

The first claim in the second assertion is proved simultaneously with the first claim in the third assertion by induction on $m$.
(We use the third assertion implicitly.)

Assume that $m\gknf A_\al(b)+l$.

Case 1. $\al=0$. Then $m<k$. Then $\pk (m-1)=m-1<m=\pk m$.

Case 2. $\al>0$. 

Case 2.1. $l>0$. $m-1\gknf A_\al(b)+l-1$.
Then for some $\beta$ we find $\pk (m-1)=\beta+l-1<\beta+l=\pk m$.

Case 2.2. $l=0$.

Case 2.2.1. $b>0$. Then for some $p$ we find $m-1\gknf A_\al(b-1)+p<m\gknf A_\al(b).$

Thus $\pk(m-1)=\psi(\om+(-2+\pk\al))+\om\cdot(\pk(b-1))+p< \psi(\om+(-2+{\pk\al)})+\om\cdot\pk(b)$.

Case 2.2.2. $b=0$.
If $\al=1$ then for some $p$ we find $\pk(m-1)=\om\cdot \pk(b-1)+p<\om\cdot \pk b=\pk m$.
Now consider the case $\al\geq 2$.
Then for some $p$ we find
$m-1=A_{\al_{k,0}}^k(0)-1=A_{\al_{k,0}}((A_{\al_{k,0}}^{k-1}(0)-1))+p$ where the latter is in $k$ normal form.
The induction hypothesis yields
$\pk (A_{\al_{k,0}}(A_{\al_{k,0}}^{k-1}(0)-1)+p )
=\psi(\om+(-2+{\pk\al_{k,0}})+\om\cdot \pk(A_{{\al_{k,0}}}^{k-1}(0)-1)+p
< \psi(\om+(-2+{\pk\al}))=\pk m$
since $\pk\al_{k,0}<\pk\al$ and $\pk\al$ is an additive principal number which is in $\psi$ normal form.
The third claim in the second assertion is proved simultaneously with the second claim in the the third assertion by a similar induction on $m$.

Assume that $m\gknf A_\al(b)+l$.

Case 1. $\al=0$. Then $m<k$. Then $\xk (m-1)=m-1<m=\xk m$.

Case 2. $\al>0$. 

Case 2.1. $l>0$. $m-1\gknf A_\al(b)+l-1$.
Then for some $\beta$ we find $\xk (m-1)=\beta+l-1<\beta+l=\xk m$.

Case 2.2. $l=0$.

Case 2.2.1. $b>0$. Then for some $p$ we find $m-1\gknf A_\al(b-1)+p<m\gknf A_\al(b).$

Thus $\xk(m-1)=\Psi(-1+{\xk\al})+\om\cdot(b-1)+p< \Psi(-1+{\xk\al})+\om \cdot b$.

Case 2.2.2. $b=0$.
If $\al=1$ then for some $p$ we find $\xk(m-1)=\om\cdot (b-1)+p<\om\cdot  b=\xk m$.
Assume now that $\al\geq 2$.

Then for some $p$ we find 
$m-1=A_{\al_{k,0}}^k(0)-1=A_{\al_{k,0}}((A_{\al_{k,0}}^{k-1}(0)-1))+p$ where the latter is in $k$ normal form.
The induction hypothesis yields
$\xk (A_{\al_{k,0}}(A_{\al_{k,0}}^{k-1}(0)-1)+p )
=\Psi(-1+{\pk\al_{k,0}})+\om\cdot (A_{\al_{k,0}})^{k-1}(0)-1)+p
< \Psi(-1+ {\pk\al})=\pk m$
since $\pk\al_{k,0}+\om^2\leq \pk\al$ because $\pk\al$ is in $\Psi$ normal form.

The second claim in the second assertion is proved by a much simpler induction on $m$.

Assume that $m\gknf A_\al(b)+l$.

Case 1. $\al=0$. Then $m<k$. Then $\ck (m-1)=m-1<m=\ck m$.

Case 2. $\al>0$. 

Case 2.1. $l>0$. $m-1\gknf A_\al(b)+l-1$.
Then for some $\beta$ we find $\ck (m-1)=\beta+l-1<\beta+l=\ck m$.

Case 2.2. $l=0$.

Case 2.2.1. $b>0$. Then for some $p$ we find $m-1\gknf A_\al(b-1)+p<m\gknf A_\al(b).$
If $\al=1$ then for some $p$ we find $\ck(m-1)=\om\cdot (1+\ck(b-1))+p<\om \cdot (1+\ck(b))=\ck m$.
If $\al>1$ then for some $\beta$ we find $\ck(m-1)=\be+\om\cdot \ck(b-1)+p< \be +\om \cdot \ck(b))=\ck m$.

Case 2.2.2. $b=0$.
Then for some $p$ we find
$m-1=A_{\al_{k,0}}^k(0)-1=A_{\al_{k,0}}((A_{\al_{k,0}}^{k-1}(0)-1))+p$ where the latter is in $k$ normal form.
If $\al=0$ or $\al=1$ we can argue as in the proof of the first claim of the second assertion.
If $\al=2$ then 
$\ck(m-1)=\ck(A_1(A_1^{k-1}(0)-1)+p)=
\om\cdot (1+\ck(A_1^{k-1}(0)-1)+p)<\om^\om=\ck(A_2(0))=\ck m$.

If $\al>2$ then
$\ck(m-1)=
\pk (A_{\al_{k,0}}(A_{\al_{k,0}}^{k-1}(0)-1)+p )=\om^{\om+(-2+\al_{k,0})}+\om\cdot \om^{\om+(-2+\al_{k,0})}+\ldots)<\om^{\om+(-2+\al)})=\ck m$
since $\al_{k,0}<\al$.

Proof of the first claim in the third assertion by induction on $m$.
Here we use the first claim of the first assertion implicitly.

If $l>0$ then the claim follows by applying the induction hypothesis to $m-1$.

Case 1. $b>0$. Then, by induction hypothesis, $n:=A_\al(b-1)$ is in $k$ normal form and so $\pk n=\psi(\om+(-2+\pk\al))+\om\cdot \pk( b-1)+p$ is in $\psi$ normal form.
Hence $G(\om+(-2+\pk\al))<\om+(-2+\pk\al)$.
Hence $\pk m=\psi(\om+(-2+\pk\al))+\om\cdot \pk b+p$ is in $\psi$ normal form.
Let us show that $\pk m$ is also in Cantor normal form.
Assume that $b$ has normal form $A_\be(c)+q$. Then $\be\leq \al$ since $m$ is in $k$-normal form.
This means that $\pk(\om+(-2+\pk \be))\leq \pk(\om+(-2+\pk \al)$. By induction hypothesis
$\pk b=\pk(\om+(-2+\pk \be))+\om \pk c+q$ is in Cantor normal form. This yields that $\pk m$  is in Cantor normal form.

Case 2. $b=0$ and $\al=\be+1$ is a successor.

Write $\be+1=\la+r$ with $\la=0$ or $\la$ a limit. Then $\pk (\be+1)=\pk \la+\pk r$ and $\pk \be=\pk \la +\pk(r-1)$.

We have to show that $G(\om+(-2+\pk(\be+1))<\om+(-2+\pk(\al)$.

We find $G(\om+(-2+\pk(\be+1))=\{1\}G\pk \la \cup G\pk r$. Since $A_\be(0)$ is in $k$ normal form the induction hypothesis yields that
$\psi(\om+(-2+\pk\be))$ is in normal form.
Hence $G\pk \la \subseteq G\pk \be\subseteq G(\om+(-2+\pk\be))<\om+(-2+\pk\be)<\om+(-2+\pk\al)$.

Now write $r=A_\ga(c)+d$ in $k$ normal form.
Then $\pk r=\psi(\om+(-2+\pk \ga))+\om\cdot \pk c+d$ is in $\psi$ normal form and in Cantor normal form by induction hypothesis.
Therefore $G\pk c\leq G \psi(\om+(-2+\pk \ga)\leq \om+(-2+\pk \ga$.
Since $m$ is in $k$ normal form we find $\ga\leq \al$. Since $\ga$ is a strict subterm of $\al$ we have $\ga< \al$
We find $G\pk r\leq \om+(-2+\pk\ga)<\om+(-2+\pk\al)$.
Hence $\pk m$ is in $\psi$ normal form.

Let us show that $\pk m$ is also in Cantor normal form.
Assume that $b$ has normal form $A_\de(e)+q$. Then $\de\leq \al$ since $m$ is in $k$-normal form.
This means that $\pk(\om+(-2+\pk \de))\leq \pk(\om+(-2+\pk \al)$. By induction hypothesis
$\pk b=\pk(\om+(-2+\pk \de))+\om \pk c+q$ is in Cantor normal form. This yields that $\pk m$  is in Cantor normal form.

Case 3. $b=0$ and $\al=\la\scbr{\om^{\be+1}\cdot r}$ is a limit where $\al[x]=\la\scbr{\om^{\be+1}\cdot (r-1)+\om^{\be}\cdot x}.$
Then $n:=A_{\al_{k-1,0}}(0)$ is in $k$ normal form and so $\pk n$ is in $\psi$ normal form by induction hypothesis.
We have $n=A_{\la[[\om^{\be+1}\cdot (r-1)+\om^\be \cdot A_{\al_{k-2,0}}(0)]]}(0)\geq \max\{mc(\al[0]),mc(\be+1),r\}$.

Let $p:=mc(\al)$ Then $p=\max\{mc(\al[0]),mc(\be+1),mc(r)\}$ and  $p\leq n$.

We find $\pk p\leq \pk n=\psi({\pk\al_{k-1,0}})$ where the letter is in $\psi$ normal form by induction hypothesis
so that $G(\psi({\pk\al_{k-1,0}}))\leq {\pk\al_{k-1,0}}$.

This yields $G(\om+(-2+{\pk\al}))\leq G(\pk p)\leq G(\pk n)=G(\psi({\pk\al_{k-1,0}}))\leq {\pk\al_{k-1,0}}< {\pk\al}$.
Hence $\pk m$ is in $\psi$ normal form.

Proof of the second claim in the third assertion by induction on $m$. The details are similar to the details in the proof of the first claim of the third assertion.

\hfill$\Box$

We define the corresponding Goodstein sequences similarly as before.
\begin{defi}
Let $m<\om$. 
\begin{enumerate}
\item
Put $m_0:=m.$
Assume recursively that $m_{l}$ is defined and $m_{l}>0$.
Then $m_{l+1}=m_l[l+3\leftarrow l+4]-1$. If $m_l=0$ then $m_{l+1}:=0$.
\item Put $\tilde{m}_0:=m.$
Assume recursively that $\tilde{m}_{l}$ is defined and $\tilde{m}_{l}>0$.
Then $\tilde{m}_{l+1}=\tilde{m}_l\scbr{l+3\leftarrow l+4}-1$. If $\tilde{m}_l=0$ then $\tilde{m}_{l+1}:=0$.
\item Put $\overline{m}_0:=m.$
Assume recursively that $\overline{m}_{l}$ is defined and $\overline{m}_{l}>0$.
Then $\overline{m}_{l+1}=\overline{m}_l\{l+3\leftarrow l+4\}-1$. If $\tilde{m}_l=0$ then $\overline{m}_{l+1}:=0$.
\end{enumerate}

\end{defi}
We are going to prove that the first principle is a giant Goodstein principle,
that the second principle is an illusionary giant Goodstein principle, and 
that somewhat surprisingly t the third principle is a giant Goodstein principle.

\begin{theo} 
\begin{enumerate}
\item For all $m<\om$ there exists an $l<\om$ such that $m_l=0.$ This is provable in $\mathrm{PRA}+\mathrm{TI}(\eta_0)$.
\item For all $m<\om$ there exists an $l<\om$ such that $\tilde{m}_l=0.$ This is provable in $\mathrm{PRA}+\mathrm{TI}(\varepsilon_0)$.
\item For all $m<\om$ there exists an $l<\om$ such that $\overline{m}_l=0.$ This is provable in $\mathrm{PRA}+\mathrm{TI}(\eta_0)$.
\end{enumerate}
\end{theo} 
Proof. Define $o(m,l):=\psi_{l+3}( m_l).$ If $m_{l+1}>0$ then by the previous lemmata
\begin{eqnarray*}
o(m,l+1)&=&\psi_{l+4}(m_{l+1})   \\
&=&\psi_{l+4} (m_l[l+3\leftarrow l+4]-1)\\
&<&\psi_{l+4} (m_l[l+3\leftarrow l+4])\\
&=&\psi_{l+3} (m_l)\\
&=&o(m,l)
\end{eqnarray*}
This proves the first assertion. The second and third assertion are proved similarly by now using $\ck$ ($\xk$ resp.) instead of $\pk$.

Let us now prove the independence results.
For this we use canonical  fundamental sequences for the elements in $OT$ which go back to Buchholz \cite{Buchholz}.
We put $0[x]:=0$ and $(\be+1)[x]:=\be$.
We put $(\psi0)[x]:=0$ and $(\psi (\be+1))[x]:=\psi\be\cdot x$.
If $\la$ is of cofinality $\om$ the we recursively put
$(\psi\la)[x]:=\psi(\la[x])$. Then $(\psi\om)[x]=\om^x$.
If $\la$ is of cofinality $\Om$ then we put
$(\psi\la)[x]:=\psi\la_{x,0}$ where $\la_{0,0}:=\la[0]$ and $\la_{l+1,0}:=\la[\psi\la_{l,0}]$.
We hereby assume that for $\al\geq\Om$ we agree by recursion on the following.
If $\al=\Om^\be\cdot\de+\ga$ and $\ga$ is a limit then the cofinality of $\al$ is the cofinality of $\ga$
and $\al[\xi]:=\Om^\be\cdot\xi+\ga[\xi]$.
If $\al=\Om^\be\cdot\de$ and $\de$ is a limit then the cofinality of $\al$ is the cofinality of $\de$
and $\al[\xi]:=\Om^\be\cdot\de[\xi]$.
If $\al=\Om^\be(\de+1)$ and $\be$ is a limit then the cofinality of $\al$ is the cofinality of $\be$
and $\al[\xi]:=\Om^\be\cdot\de+\Om^{\be[\xi]}$.
If $\al=\Om^{\be+1}(\de+1)$ then $\al$ has cofinality $\Om$ and 
$\al[\xi]=\Om^{\be+1}\de+\Om^{\be}\cdot \xi.$

These fundamental sequences have the Bachmann property \cite{Weiermanna}. Moreover, by \cite{Weiermanna} we have that 
if $\psi\al$ is in $\psi$ normal form then $(\psi\al)[x]$ is in $\psi$ normal form, too. These two results go back to Buchholz \cite{Buchholz}.

The fundamental sequences for $OT'$ are defined in complete analogy. The only difference is the clause $\Psi(\al+1)[x]=\Psi\al+\om\cdot x$.
As before one can prove that if $\Psi\al$ is in $\Psi$ normal form then $(\Psi\al)[x]$ is in $\Psi$ normal form, too.

We call a natural number $m$ of $k$-successor type if $m\gknf A_\al(b)+l$ where $\al=0$, or $\al>0$ and $l>0$.
For those $m$ the ordinal $\pk m$ is a successor ordinal.
We call a natural number $m$ of $k$-limit type if $m\gknf A_\al(b)$ with $\al\geq 1$.
For those $m$ the ordinal $\pk m$ is a limit ordinal of countable cofinality.

For technical reasons we need a specific description of ordinals below $\eo$ in terms of certain place holders.

\begin{defi}\label{defelf}
We define by recursion on $\al$ a context $\la_k(\al)$ for $0<\al<\eo$ which are not of the form $\be+q$ with $\be\in Lim\cup\{0\}$ and $q$ is of $k$ successor type.
Assume that 
$\al=\om^{\al_1}\cdot m_1+\cdots+\om^{\al_n}\cdot m_n$ where $\al_1>\ldots>\al_n$ and $0<m_1,\ldots,m_n$.

Case 1. $m_n$ is of $k$ limit type. Then
$\la_k(\al):=\om^{\al_1}\cdot m_1+\cdots+\om^{\al_{n-1}}\cdot m_{n-1}+\scbr{\cdot}$.
Then $\la_k(\al)\scbr{\om^{\al_n}\cdot m_n}=\al$.

Case 2. $m_n$ is of $k$ successor type. We have excluded the case $\al_n=0$ by assumption.

Case 2.1.  $\al_n=\be+q$ with $\be\in Lim\cup\{0\}$ and $q$ is of $k$ successor type.

Then $\la_k(\al):=\om^{\al_1}\cdot m_1+\cdots+\om^{\al_{n-1}}\cdot m_{n-1}+\scbr{\cdot}$.
Then $\la_k(\al)\scbr{\om^{\al_n}\cdot m_n}=\al$.

Case 2.2.  $\al_n\in Lim$ or $\al_n=\be+q$ with $\be\in Lim\cup\{0\}$ and $q$ is of $k$ limit type.
By recursion we can assume that $\la_k(\al_n)$ is defined
and that $\la_k(\al)\scbr{\om^\ga\cdot p}=\al$ where $p$ is of $k$ limit type, or $p$  is of $k$ successor type and 
$\ga=\be+r$ with $\be\in Lim\cup\{0\}$ and $r$ is of $k$ successor type. 

$\la_k(\al):=\om^{\al_1}\cdot m_1+\cdots+\om^{\al_{n}}\cdot (m_{n}-1)+\om^{\la_k(\al_n)}$.
Then $\la_k(\al)\scbr{\om^\ga\cdot p}=\al$.
\end{defi}
Then $\la_k(\al)\scbr{\om^\ga\cdot p}=\al$
where $p$ is of $k$ limit type, or $p$ is of $k$ successor type and 
$\ga=\be+q$ with $\be\in Lim\cup\{0\}$ and $q$ is of $k$ successor type. 
It is easy to show that $(\la_k(\al)\scbr{\om^\ga\cdot p})'=(\la_k(\al))'\scbr{\om^{\ga'}\cdot p'}$
and $\pk(\la_k(\al)\scbr{\om^\ga\cdot p})=(\pk(\la_k))\scbr{\Om^{\pk\ga}\cdot \pk p}$.
This is because if $m$ is of $k$ successor type then $(m-1)'+1=m'$ and $\pk(m-1)+1=\pk m$.

If $\la_k(\al)\scbr{\om^\ga\cdot p}=\al$ and $p$ is of $k$ limit type then $\pk\al$ is a limit of countable cofinality
and $(\pk\al)[x]=(\pk(\la_k))\scbr{\Om^{\pk\ga}\cdot (\pk p)[x]}$.

\begin{lem}\label{lem15}. Assume that $m>0$.
\begin{enumerate} 
\item $\pk m>\pskpe (m'-1)\geq (\pk m)[k-2]$,
\item $\ck m>\cskpe (m'-1)\geq (\ck m)[k-2]$,
\item $\xk m >\xskpe (m'-1)\geq (\xk m)[k-1]$.
\end{enumerate}
\end{lem}
Proof. Clearly $\pk m=\pskpe (m')>\pskpe (m'-1).$  This argument works also for $\ck$ and $\xk$
so that we only need to show the second inequality in every assertion.

So let us first prove the second inequality of the first assertion.
If $0<m<k$ then $\pskpe(m'-1)=m-1= (\pk m)[k-2].$

Assume that $m\gknf A_\al(b)+l\geq k$. Then $\al>0$.

Case 1. $l>0$. Then $m-1\gknf A_\al(b)+l-1$.
Then $(\pk m)[k-2]=(\pk (A_\al(b))+l)[k-2]=\pk (A_\al(b))+l-1$ and 
$\pskpe (m'-1)=\pskpe(B_{\al'}(b')+l-1)=\pk (A_\al(b))+l-1$.

Case 2. $l=0$.

Let us first consider the case $\al=1$.
If $b=0$ then 
The $\pk A_1(0)[k-2]=\om[k-2]\leq k=\pskpe (k+1-1)=\pskpe (m'-1)$.

If $b>0$ then
$\pk A_1(b)[k-2]=(\om(1+\pk b))[k-2]\leq \om(1+\pk b[k])+k$.
Moreover for some $p\geq k$ the induction hypothesis yields $\pskpe (B_1(b')-1)\geq \pskpe (B_1(b'-1)+p)\geq \om (1+\pskpe(b'-1))+p\geq \om(1+\pk b[k-2])+k$.

Let us second consider the case $\al=2$.
If $b=0$ then 
The $\pk A_2(0)[k]=\psi(\om+(-2+\pk 2))[k-2]=(\psi\om)[k-2]=\psi(\om[k-2])\leq \om^k
\leq \pskpe(B_1^k(0))\leq \pskpe(B_1^{k+1}(0)-1)=\pskpe (m'-1)$.

If $b>0$ then
$\pk A_2(b)[k-2]=(\psi\om+ \om(\pk b))[k-2]\leq \psi\om+\om(1+\pk b[k-2])+k$.
Moreover for some $p\geq k$ the induction hypothesis yields $\pskpe B_2(b')-1\geq \pskpe (B_2(b'-1)+p)\geq \psi \om+\om (\pskpe(b'-1))+p\geq \psi\om+\om(1+\pk b[k-2])+k$.

From now on we assume that $\al\geq 3$.

Case 2.1. $b>0$. 
Then for some $l\geq k$ $m'-1=B_{\al'}(b')-1=B_{\al'}(b'-1)+l$ where the latter is in $k+1$ normal form.

Here the induction hypothesis yields
$(\pk m)[k-2]=(\psi(\om+(-2+{\pk \al}))+\om \cdot(\pk b))[k-2] \leq \psi(-2+{\pk \al}+\om \cdot(\pk b[k-2]))+k
\leq \psi(\om+(-\pskpe 2+{\pskpe \al'}))+\om \cdot(\pskpe(b'-1))+p=
 \pskpe(B_{\al'}(b'-1)+p)
= \pskpe(B_{\al'}(b')-1)=\pskpe (m'-1)$.

Case 2.2. $b=0$.

Then $B_{\al'}(0)-1=B_{\al'_{k+1,0}}(0)^{k+1}-1>B^{k}_{\al'_{k,0}}(0)$ where the latter is in $k+1$ normal form.

Case 2.2.1. $\al$ is a limit. 

We now analyse the specific forms for $\al$ according to Definition \ref{defelf}.

Case 2.2.1.1. $\al=\la_k(\al)\scbr{\om^{\be+p}\cdot q}$ where $p$ and $q$ are of $k$ successor type and $\be$ is a limit or zero.
Then $\pk\al=(\pk\la_k(\al))\scbr{\Om^{\pk\be+\pk p}\cdot \pk q}$ has uncountable cofinality since $\pk p=\pk (p-1)+1$ and $\pk q=\pk (q-1)+1$.

Then $(\pk m)[k]=(\psi(\om+(-2+{\pk \al})))[k] =(\psi({\pk \al}))[k]=  \psi(({\pk(\al)})_{k,0})$ and
$\pskpe (m'-1)=\pskpe(B_{\al'}(0)-1)=\pskpe(B^{k+1}_{(\al')_{k+1,0}}(0)-1)\geq \pskpe(B_{(\al')_{k,0}}(0))=\psi(\om+(-2+{\pskpe(\al'_{k,0}})))=\psi({\pskpe(\al'_{k,0})})$ since
$B_{(\al')_{k,0}}(0)$ is in $k+1$ normal form because $B_{\al'}(0)$ is in $k+1$ normal form.

We claim that 
$({\pk\al})_{l,0})={\pskpe(\al'_{l,k+1,0}})$ for $l\leq k$. This then yields the assertion in this case.

We prove this by induction on $l$. 
Note first that $\al'=(\la_k(\al))'\scbr{\om^{\be'+p'}\cdot q'}$ where $p'=(p-1)'+1$ and $q'=(q-1)'+1$.
Further note that $\pk\al=(\pk\la_k(\al))\scbr{\Om^{\pk\be+\pk p}\cdot \pk q}$. 

Assume that $l=0$.
We find $\al'_{0,0}=\al'[0]=(\la_k(\al)')\scbr{\om^{\be'+p'}\cdot (q-1)'}$
and $(\pk\al)[0]=(\pk\la_k(\al))\scbr{\Om^{\pk\be+p}\cdot \pk(q-1)}$.
This yields 
\begin{eqnarray*}
&&(\pk\al)_{0,0}\\
&=&(\pk\al)[0]\\
&=&(\pk\la_k(\al))\scbr{\Om^{\pk\be+\pk p}\cdot \pk (q-1)}\\
&=&{(\pskpe(\la_k(\al))')\scbr{\Om^{\pskpe\be'+\pskpe p'}\cdot \pskpe(q-1)'}}\\
&=&\pskpe(\al'[0])\\
&=&{\pskpe(\al'_{0,k+1,0}}).
\end{eqnarray*}
Now assume assume that the claim is true for $l$.
Then
\begin{eqnarray*}
&&({\pk\al})_{l+1,0}\\
&=&{\pk\la_k(\al))\scbr{\Om^{\pk\be+\pk p}\cdot \pk (q-1)+\Om^{\pk\be+\pk(p-1)}\cdot  \psi(({\pk\al})_{l,0} )}}\\
&=&{(\pskpe(\la_k(\al))')\scbr{\Om^{\pskpe\be'+\pskpe p'}\cdot \pskpe(q-1)'+\Om^{\pskpe\be'+\pskpe(p-1)'}\cdot  \psi{(\pskpe( \al'_{l,k+1,0}}) )}}\\
&=&{\pskpe(\la_k(\al)'\scbr{\om^{\be'+\pskpe p'}\cdot \pskpe(q-1)'+   \om^{\be'+\pskpe(p-1)'}\cdot B_{ \al'_{l,k+1,0} }(0)}}\\\
&=&{\pskpe(\al'_{l+1,k+1,0})}.
\end{eqnarray*}

Case 2.2.1.2. $\al=\la_k(\al)\scbr{\om^\ga\cdot q}$ with $q$ of $ k$ limit type. Then $q'$ of $ k+1$ limit type
Then $\pk\al=(\pk\la_k(\al))\scbr{\Om^{\pk{\ga}}\cdot \pk q}$ is a limit of countable cofinality
and $\al'=(\la_k(\al))'\scbr{ \om^{\ga'}\cdot q'}$ where $q'\geq k+1$.

Since $\al$ is a limity we have $(-2+\pk \al)[k-2]=-2+\pk \al[k-2]$ and we find
\begin{eqnarray*}
&&(\pk m)[k]\\
&=&(\psi{(\om+(-2+\pk\al}))[k]\\
&=&\psi{(\om+(-2 +(\pk\al))[k])}\\
&=&\psi{(\om+(-2+\pk\la_k(\al)\scbr{\Om^{\pk\ga}\cdot (\pk q))[k-2]}})\\
&\leq &\psi(\om+(-2+{\pskpe(\la_k(\al)'\scbr{\Om^{\pskpe\ga'}\cdot \pskpe (q'-1)}}))))\\
&= &\pskpe(B_{\la_k(\al)'\scbr{\om^{\ga'}\cdot  (q'-1)}}(0))\\
&\leq &\pskpe(B_{\al'}(0)-1)
\end{eqnarray*}

since $\la_k(\al)'\scbr{\om^{\ga'}\cdot  (q'-1)}\preceq_1 \al'$.

Case 2.2.2. $\al$ is a successor say $\al=\be+r$ with $r>0$ and $\be$ is zero or a limit. If $\be=0$ then $r\geq 3$.

If $\al=k$ then $(\pk m)[k-2]=\pk(\om+(-2+\pk k))[k-2]=\psi(\om+\om)[k-2]=\psi(\om+(-2+k))=\psi(\om+(-2+\pskpe(\al'-1)))=\pskpe B_{\al'-1}(0)\leq \pskpe (B_{\al'}(0)-1).$
Now assume $\al\not= k$ so that $(-2+\pk \al)[k-2]=-2+\pk \al[k-2]$.

Here the induction hypothesis yields

\begin{eqnarray*}
&&(\pk m)[k-2]\\
&=&\psi(\om+(-\pk2+{\pk \al}))[k-2]\\
&\leq &
\psi{\om+(-2+\pk(\be)+\pk(r)[k-2]))}\cdot k \\
&\leq &\psi(\om+(-2 +{\pskpe(\be')+\pskpe(r'-1)})\cdot k\\
&=&\psi(\om+(-2+{\pskpe(\be'+r'-1)))}\cdot k \\
&\leq &\pskpe(B_{\be'+r'-1}^k(0))\\
&\leq &\pskpe(B_{\al'_{k+1,0}}^{k+1}(0)-1)\\
&<& \pskpe(B_{\al'}(0)-1)=\pskpe (m'-1).
\end{eqnarray*}

The second equality in the third assertion is proved by a similar induction on $m$.
The differences will be very small and so we skip most of the proof.
Let $m\gknf A_\al(k,b)+l$. The case $l>0$ is as before. The case $b>0$ is also similar as before. 
But for $b>0$ we  use the fact that $\Psi(\ga+1)=(\Psi\ga)+ \om^2$ to model the $\xk$ interpretation.
The fact that $\Psi(\ga+1)=(\Psi\ga)+ \om^2$ is also used to model the $\xk$ interpretation in the case $\al=\be+r$ with $r>0$ and $\be$ is zero or a limit. 
The fact that we use $k-1$ as argument of the fundamental sequence in the assertion has to do with the degenerate case $\al=k,b=0,l=0$ where we need
 $(\om^2\cdot (-1+\om))[k-1]=\om^2\cdot (k-1)=\om^2\cdot (-1+\pskpe(k'))$.

Let us now prove second equality in the second assertion.

If $0<m<k$ then $\cskpe(m''-1)=m-1= (\ck m)[k-2].$

Assume that $m\gknf A_\al(k,b)+l\geq k$. Then $\al>0$.

Case 1. $l>0$. Then $m-1\gknf A_\al(k,b)+l-1$.
Then $(\ck m)[k-2]=(\ck (A_\al(k,b))+l)[k-2]=\ck (A_\al(k,b))+l-1$ and 
$\cskpe (m''-1)=\cskpe(B_{\al}(b'')+l-1)=\ck (A_\al(k,b))+l-1$.

Case 2. $l=0$.

Case 2.1. $b>0$. 

If $\al=1$ then the induction hypothesis yields $(\ck m)[k-2]=\om(1+\ck b)[k-2]\leq \om(1+(\ck b)[k-2])+k\leq \om(1+(\cskpe (b''-1)))+k= \cskpe(B_{\al'}(b''-1)+k\leq \cskpe(B_{\al}(b'')-1)=\cskpe (m''-1)$.

If $\al>1$ then the induction hypothesis yields
$(\ck m)[k-2]=(\om^{\om+(-2+ \al)} +\om\cdot \ck b)[k-2] \leq   \om^{\om+(-2+\al)} +\om\cdot ((\ck b)[k-2] )+k\leq 
\om^{\om+(-2+ \al)} +\om\cdot ((\cskpe b''-1))+k=\cskpe(B_{\al}(b''-1)+k\leq \cskpe(B_{\al}(b'')-1)=\cskpe (m''-1)$.

Case 2.2. $b=0$.
If $\al=1$ then $m=A_1(0)=k$ and $(\ck m)[k-2]=\om[k-2]\leq k= \cskpe(k+1-1)=(\cskpe (m''-1)$.

If $\al=2$ then $(\ck m)[k]=(\om^{\om+(-2+2})+\om\cdot 0)[k-2]=\om^\om[k-2]\leq \om^k=\om(1+\cskpe(B_1^{k-1}(0))=\pskpe(B_1^k(0))\leq 
\cskpe (B_1^{k+1}(0)-1)=\cskpe (m''-1)$. Note that $\cskpe B_1^l(0)=\om^l$ holds by induction on $l$ for $0<l\leq k$.

If $\al>2$ then 

\begin{eqnarray*}
&&(\ck m)[k]=(\om^{\om+(-2+ \al)})[k-2]\\
&\leq& (\om^{\om+-2+(\al)[k])})\cdot k\\
&\leq& (\om^{\om+(-2+\al_{k,0})})\cdot k\\
&=&\cskpe(B_{\al_{k,0}}^k(0))\\
&\leq &\cskpe(B_{\al}(0)-1)\\
&\leq &\cskpe(B_{\al}(0)''-1)\\
&=&\cskpe (m''-1)
\end{eqnarray*}

Note that $\pskpe (B_{\al_{k,0}}^l(0))=\om^{\om+(-2+\al_{k,0})}\cdot l$ holds by induction on $l$ for $0<l\leq k$.

\begin{theo}

\begin{enumerate}
\item $\mathrm{ID}_1\not \vdash (\forall m) (\exists l )[m_l=0]$.
\item $\mathrm{PA}\not \vdash (\forall m )(\exists l )[\tilde{m}_l=0]$.
\item $\mathrm{PA}\not \vdash (\forall m) (\exists l )[\overline{m}_l=0]$.
\end{enumerate}
\end{theo}

Proof. Let $\om_1:=\om$ and $\om_{r+1}:=\om^{\om_r}$. 
Let $m(r):=A_{\om_r}(3,0).$
Let $\Om_1:=\Om$ and $\Om_{r+1}:=\Om^{\Om_r}$.

Then $\psi_3(m(r))=\psi(\Om_r)$ for $r\geq 1$.
We claim that $o(m(k),l)\geq_1\psi(\Om_r)[1]\ldots [l]$.
Proof of the claim. Write $m$ for $m(r)$.
For $o(m,l)>0$ we have $o(m,l)>o(m,l+1)=\psi_{l+4}(m_l[l+3\leftarrow l+4]-1)\geq (\psi_{l+3}(m_l))[l+1]=(o(m,l)[l+1]$.
The Bachmann property yields $o(m,l+1)\geq _1 (o(m,l)[l+1]$.
The induction hypothesis yields $o(m,l)\geq_1 \psi(\Om_r)[1]\ldots [l]$
hence $o(m,l)[l+1]\geq_1 \psi(\Om_r)[2]\ldots [l][l+1]$.
Therefore $o(m,l+1)\geq_1 (o(m,l)[l+1]\geq_1  \psi(\Om_r)[1]\ldots [l][l+1]$.

Therefore the least $l$ such that $o(m,l)=0$ is at least as big as the least $l$ such that 
$\psi(\Om_r)[1]\ldots [l]=0$. The result follows from $ID_1\not\vdash \forall r \exists l (\psi(\Om_r))[2]\ldots [l+1]=0$.
(See, for example, \cite{Buchholzb} for a proof.)

The second assertion follows similarly. 
We see $\chi_3(m(r))=\om^{\om_r}$ for $r\geq 2$.
The result follows from $\mathrm{PA }\not\vdash (\forall r)( \exists l )[({\om_{r+1}})[1]\ldots [l]=0]$.

The third assertion follows similarly by using $m(r)$ again and $\xi_3$.\hfill$\Box$

So the principles $ (\forall m) (\exists l )[m_l=0]$ and $(\forall m) (\exists l )[\overline{m}_l=0]$ turn out to be giant Goodstein principles and the principle
$(\forall m )(\exists l )[\tilde{m}_l=0]$  turns out to be a (non trivial) illusionary giant Goodstein principle.

A weaker Goodstein principle can also be obtained by performing a trivial base change in the iteration parameter $k$.
For $m\gknf A_\al(k,b)+l\geq k$ let $m'''':=A_\al(k+1,b)+l$. A Goodstein principle base on this definition becomes 
provable in a theory in which we can define this Goodstein process and in which we can work with the ordinal assignment $ord(m):=\om^\al+\om \cdot b+l$.


\begin{thebibliography}{1}

\bibitem{Arai} T.~Arai, T, D.~Fernandez-Duque, S.~Wainer, A.~Weiermann: Predicatively unprovable
termination of the Ackermannian Goodstein principle. Proceedings of the
American Mathematical Society 148 (2020), no. 8, 3567--3582. https://doi.org/10.1090/proc/14813 
\bibitem{Buchholza} W.~Buchholz, S.~Feferman, W.~Pohlers, W.~Sieg: Iterated inductive definitions and subsystems of analysis: recent proof-theoretical studies. Lecture Notes in Mathematics, 897. Springer-Verlag, Berlin-New York, 1981. v+383 pp. https://10.1007/BFb0091894
\bibitem{Buchholz} W.~Buchholz. 
A new system of proof-theoretic ordinal functions.
Ann. Pure Appl. Logic 32 (1986), no. 3, 195--207. https://doi.org/10.1016/0168-0072(86)90052-7
\bibitem{Buchholzb} W.~Buchholz. An independence result for $(\Pi^1_1-\mathrm{CA})+\mathrm{BI}$. Annals of  Pure and Applied Logic 33 (1987), no. 2, 131--155.
https://doi.org/10.1016/0168-0072(87)90078-9
\bibitem{BCW} W.~Buchholz. A.~ Cichon, A.~Weiermann:
A uniform approach to fundamental sequences and hierarchies. 
Math. Logic Quart. 40 (1994), no. 2, 273--286. https://doi.org/10.1002/malq.19940400212
\bibitem{Cichon} E.A.~Cichon.
{A short proof of two recently discovered independence results
using recursion theoretic methods.}
Proceedings of the AMS 87 (1983) 704-706. https://doi.org/10.1090/S0002-9939-1983-0687646-0 
\bibitem{Fernandez} D.~Fernandez-Duque, A.~Weiermann:  Ackermannian Goodstein Sequences of Intermediate Growth. In: Anselmo M., Della Vedova G., Manea F., Pauly A. (eds) Beyond the Horizon of Computability. CiE 2020. Lecture Notes in Computer Science, vol 12098. Springer, Cham. 163-174 https://doi.org/10.1007/978-3-030-51466-2-14 
\bibitem{walk} D.~Fernandez-Duque, A.~Weiermann: A walk with Goodstein. (submitted) https://arxiv.org/abs/2004.09110
\bibitem{Goodsteina} R.L.~Goodstein.
On the restricted ordinal theorem.
J. Symbolic Logic 9, (1944). 33--41. https://www.jstor.org/stable/2268019 
\bibitem{Goodsteinb} 
R.L.~Goodstein.
Transfinite ordinals in recursive number theory.
J. Symbolic Logic 12, (1947). 123--129. https://www.jstor.org/stable/2266486 

\bibitem{Kirby}
L.~Kirby and J.~Paris.
Accessible independence results for Peano arithmetic.
Bull. London Math. Soc. 14 (1982), no. 4, 285--293.  https://doi.org/10.1112/blms/14.4.285
\bibitem{Weiermannb} A.~Weiermann. Classifying the provably total functions of PA. Bull. Symbolic Logic 12 (2006), no. 2, 177--190. https://www.jstor.org/stable/4617258
\bibitem{Weiermanna} A.~Weiermann. Investigations on slow versus fast growing: how to majorize slow growing functions nontrivially by fast growing ones. 
Arch. Math. Logic 34 (1995), no. 5, 313--330. https://doi.org/10.1007/BF01387511
\bibitem{Weiermann} A.~Weiermann. Ackermannian Goodstein principles for first order Peano arithmetic.  Sets and computations, 157--181,
Lect. Notes Ser. Inst. Math. Sci. Natl. Univ. Singap., 33, World Sci. Publ., Hackensack, NJ, 2018. https://doi.org/10.1142/9789813223523-0007

\end{thebibliography}
\end{document}